\documentclass[10pt]{article}
\usepackage{amsmath,amssymb,amsthm}
\usepackage{hyperref}
\usepackage{geometry}
\usepackage{graphicx}
\newcommand{\ee}{\mathrm{e}}
\newcommand{\ii}{\mathrm{i}}
\newcommand{\eps}{\varepsilon}
\newcommand{\la}{\lambda}
\newcommand{\hLa}{\widehat{\Lambda}}
\newcommand{\BC}{\mathbb{C}}
\newcommand{\BE}{\mathbb{E}}
\newcommand{\BN}{\mathbb{N}}
\newcommand{\BP}{\mathbb{P}}
\newcommand{\BR}{\mathbb{R}}
\newcommand{\BZ}{\mathbb{Z}}
\newcommand{\CC}{\mathcal{C}}

\newcommand{\CP}{\mathcal{P}}
\newcommand{\CB}{\mathcal{B}}
\newcommand{\hp}{\widehat{p}}
\newcommand{\hq}{\widehat{q}}
\newcommand{\hN}{\widehat{N}}
\newcommand{\dd}{\mathrm{d}}
\newcommand{\ta}{\widetilde{a}}
\newcommand{\tb}{\widetilde{b}}
\newcommand{\tK}{\widetilde{K}}
\newcommand{\tR}{\widetilde{R}}
\newcommand{\tla}{\widetilde{\lambda}}
\newcommand{\tPi}{\widetilde{\Pi}}

\newcommand{\trho}{\widetilde{\rho}}
\newcommand{\tp}{{\bar p}}
\newcommand{\tq}{{\bar q}}
\newcommand{\sbar}{{\bar s}}
\newcommand{\boldone}{\boldsymbol{1}}

\newcommand{\eqdef}{\mathrel{{:}{=}}}

\newcommand{\mgc}[1]{M^{\text{gc}}_{#1}}

\newtheorem{thm}{Theorem}
\newtheorem{proposition}[thm]{Proposition}
\newtheorem{lemma}[thm]{Lemma}
\newtheorem{corollary}[thm]{Corollary}
\theoremstyle{definition}
\newtheorem*{remark*}{Remark}
\newtheorem{definition}[thm]{Definition}
\numberwithin{thm}{section}
\numberwithin{equation}{section}

\title{Growing integer partitions with uniform marginals\\ and the equivalence of partition ensembles}
\author{Yuri Yakubovich\footnote{St.~Petersburg State University, 7/9 Universitetskaya nab., St.~Petersburg, 199034, Russia.\\ E-mail: \texttt{y.yakubovich@spbu.ru}}}
\date{\hfill{\small\textit{To the memory of Anatoly Vershik}}}

\begin{document}
\maketitle
\begin{abstract}
We present an explicit construction of a Markovian random growth process on integer partitions such that 
given it visits some level $n$, it passes through any partition $\lambda$ of $n$ with equal 
probabilities.  The construction has continuous time, but we also investigate its discrete time jump chain.
The jump probabilities are given by explicit but complicated expressions, so we 
find their asymptotic behavior as the partition becomes large.  This allows us to explain how the limit shape is formed. 

Using the known connection
of the considered probabilistic objects to Poisson point processes, we give an alternative description 
of the partition growth process in these terms.  Then we apply the constructed growth process to 
find sufficient conditions for a phenomenon known as equivalence of two ensembles of random partitions
for a finite number of partition characteristics. This result allows to show
that counts of odd and even parts in a random partition of $n$ are asymptotically
independent as $n\to\infty$ and to find their limiting distributions, which
are, somewhat surprisingly, different.
\end{abstract}


\section{Introduction}

A partition of a non-negative integer $n$ is an infinite non-increasing sequence 
$\la=(\la_1,\la_2,\dots)$ of non-negative integers, usually called parts, that sum up to $n$.  Since the sum is
finite, the sequence is eventually zero, and the minimal number $\ell=\ell(\la)$ such that $\la_{\ell+1}=\la_{\ell+2}=\dots=0$
is called the \textit{length} of partition.  It is also common to identify $\la$ with the finite sequence $(\la_1,\dots,\la_\ell)$ in this case.   We
write $N(\la)=\sum_j\la_j$ for the \textit{weight} of a partition $\la$, that is the sum of its parts, and also write 
$\la\vdash n$\label{vdash} to denote that $\la$ is a partition of $n$.  We denote $\CP_n$ the set of partitions of $n$ and $\CP=\cup_{n\ge 0}\CP_n$\label{CP} the
set of all integer partitions, including the empty partition $\varnothing$\label{empty-partition} of $0$.  Sometimes it is convenient to parameterize partitions with part counts $C_k(\la)=\#\{j:\la_j=k\}$, $k=1,2,\dots$. 
Then 
\begin{equation}\label{N-via-Rk}
\ell(\la)=\sum_{k=1}^\infty C_k(\la),\qquad N(\la)=\sum_{k=1}^\infty kC_k(\la).
\end{equation} 

An integer partition $\la$ can be visualized by the so called Young diagram  $Y_\la$: it is a subset of the positive
quadrant $\BR_+^2$ comprised of unit squares, called \textit{boxes}, organized in rows of lengths $\la_1,\dots,\la_\ell$
and placed one above another left-aligned to the ordinate axis.  
Slightly abusing the notation, we identify the Young diagram $Y_\lambda$ with a function of $x\ge 0$ defined as
\begin{equation}\label{def-Y}
Y_\lambda(x)=\sum_{k\ge x}C_k(\lambda),\quad (x>0);\qquad Y_\lambda(0)=\lim_{x\downarrow0}Y_\lambda(x)=\ell(\la),
\end{equation}
so that the Young diagram is the subgraph of this function.

The inclusion partial order on Young diagrams is inherited
by partitions, and we write $\mu\nearrow\la$ (or $\la\nwarrow\mu$)\label{nearrow} to denote that $Y_\mu\subset Y_\la$ and $ Y_\la\setminus Y_\mu$ 
is just one box.  In other words, $\la$ covers $\mu$ in the inherited partial order.  Naturally, it is possible only if $N(\la)=N(\mu)+1$,
and $N$ is a rank function in this poset.
Its Hasse diagram is known as the \textit{Young graph} which we also denote by $\CP$.  We refer to $\CP_n\subset \CP$ as the $n$-th 
\textit{level} of the Young graph, that is the set of all vertices of rank $n$. 

Partitions appear naturally in many contexts, usually as block sizes of decomposable combinatorial structures. 
For instance, permutations of $n$ objects can be decomposed into cycles, and cycle lengths form a partition of $n$.
A uniform probability measure on permutations of $n$ objects then is projected to a probability measure $M_n^{\text{Ewens(1)}}$
on partitions of~$n$, known as the Ewens measure (with parameter 1) \cite{CSP}.  
Partitions of $n$ are also known to index irreducible representations of the symmetric group 
on $n$ objects (see, e.g., \cite{James}), and the Plancherel measure \cite{Vershik-Kerov77} is $M^{\text{Plancherel}}_n(\la)=(\dim \la)^2/n!$, where
$\dim \la$ is the dimension of the irreducible representation indexed by $\la$.

In each of these two examples a probability measure $M_n$ on $\CP_n$ is defined for each $n=0,1,2,\dots$. 
This sequence of measures can be also obtained by a Markovian stochastic procedure sending a partition
$\mu\vdash n$ to a partition $\la\vdash n+1$ such that $\la\nwarrow\mu$
which sends the measure $M_n$ to $M_{n+1}$.   This procedure is particularly easy for $M_n^{\text{Ewens(1)}}$: 
given a partition $\mu\vdash n$ we either increase a part $\mu_j$, $j=1,\dots,\ell(\mu)$, by 1 with
probability $\mu_j/(n+1)$, or add a new part 1 with probability $1/(n+1)$ \cite{CSP}. If $\mu$ is taken at random
with the distribution $M_n^{\text{Ewens(1)}}$ then a new partition $\la$ has the distribution $M_{n+1}^{\text{Ewens(1)}}$.
The procedure for $M_n^{\text{Plancherel}}$ is also very explicit but harder to describe, see \cite{Kerov93} and \cite{Olsh-Ivanov}.  

Many other examples of growth procedures for various measures on the Young graph or on more general graphs are also known.
For instance, the procedure leading to the Ewens measure with parameter~1 described above can be generalized much 
further to deal with measures on integer partitions which are projections of exchangeable partitions of finite sets,
including the most important example of a two parametric extension of the Ewens measure, see~\cite{CSP}.
An opposite approach, when some natural procedure of diagram growth is introduced and then the resulting
measures on partitions are investigated, has been also employed by some researchers, see, e.g.,\ \cite{Krapivsky,Eriksson,Fatkullin}.

\smallskip

However for another natural example, that is the uniform measures $M_n^{\text{uni}}$\label{Mnuni} on the sets $\CP_n$ of partitions of $n$, 
no such procedure is described in the literature, to the best of author's knowledge.  Moreover, it is not known if such procedure exists,
see Section \ref{sec:adding-a-cell}.  In this paper we propose another stochastic Markovian procedure of Young diagram growth
which leads to the uniform distribution.
While in the previous examples just one box is added to the Young diagram on each step, in our procedure a rectangular 
block of boxes is added.  So the growing diagram jumps over some levels of the Young graph $\CP$. However given
the growing diagram visits the $n$-th level $\CP_n$, the conditional probability that it passes through a diagram $Y_\la$, $\la\vdash n$,
does not depend on $\la$, hence it is the uniform measure $M_n^{\text{uni}}$.

Our procedure thus can not be treated as an increasing random walk on the Young graph.  As any Markov chain, 
our chain still can be understood as a random walk on a different graph, but this graph is not graded, and each its
vertex has infinitely many neighbors, although intervals are finite. 

The construction we introduce is based on a new look on the so-called Fristedt's conditional device \cite{Fristedt},
which basically allows to describe a uniform random partition of $n$ by a sequence of independent random variables
conditioned that their certain linear combination equals $n$ (or, in other words, that the random partition lies in $\CP_n$), 
see details in Section~\ref{sec:gce}.  The distributions
of these independent random variables depend on a parameter $q\in(0,1)$, and a proper choice $q=q_n$
of this parameter allows to mimic the behavior of the uniform probability measure on $\CP_n$.  On the other
hand, independence makes it easy to establish limiting results for various functionals of partitions.  Here
two alternatives can arise: either the conditioning on lying in $\CP_n$ does not affect the limiting
distribution and it suffices just to plug the correct value $q_n$ for the parameter to get the corresponding result 
for $M_n^{\text{uni}}$, or the limiting distribution is changed by the conditioning and additional
efforts are needed to translate the result into a conditional claim.  In the former case we say that the
\textit{unconditional and conditional measures are equivalent} for the functional considered.  A typical example
of a functional for which measures are equivalent is the limit shape (where the limiting distribution is 
degenerate). Another example is the length of random partition, which after an appropriate scaling approaches Gumbel's distribution.
More fine functionals like fluctuations around the limit shape or the size of the Durfee square are known \cite{Mutafchiev} to behave
differently before and after conditioning.

To illustrate a utility of the constructed partition growth process we use it to provide a sufficient condition on a functional
which guarantees that the unconditional and conditional measures are equivalent for this functional (Theorem~\ref{th:equivalence}).  
Looking forward to future applications,
we formulate it for joint distributions of several functionals.
In this paper we use Theorem~\ref{th:equivalence} to
find a joint limit distribution for separate counts of odd and even parts
in the uniform random partition of large integer (Theorem~\ref{th:odd-even}).

\smallskip 

The rest of the paper is organized as follows. In Section~\ref{sec:gce} we review results by Fristedt \cite{Fristedt}
which provide a base for our construction, and then introduce a basic continuous time decreasing Markov process $\la(t)$ on partitions.
In Section~\ref{sec:jumps} a jump chain of $\la(t)$ and its time-reversed chain are investigated. The forward (decreasing) jump chain has 
a simple combinatorial description (Proposition \ref{prop:jump-forward}), while the 
backward (increasing) jump chain provides a Markovian procedure of partition growth (see Theorem \ref{th:hat-Lambda}), 
which is the primary subject of this paper.
Section~\ref{sec:Ppp} establishes a connection of this model with the Poisson point processes, which is somewhat simpler
for explicit calculations.  In Section~\ref{sec:ls} we explain why the limit shape is formed for random partitions of 
a large integer.  Section~\ref{sec:adding-a-cell} contains some notes about the possibility to construct an alternative
growth model by adding one box to a Young diagram on each step, which could be easier to use for understanding the
properties of growing partitions; however, our results here are rather limited.  Our principal result about
equivalence of ensembles for a functional on partitions (Theorem~\ref{th:equivalence}) is formulated and proved
in Section~\ref{sec:equiv}.  Section~\ref{sec:oddeven} presents an application 
of this technique: we study a joint
distribution of odd and even part counts in the uniform random partition
of a large number.  The main part of the paper is concluded with some remarks in 
Section~\ref{sec:remarks}.

The paper contains also two appendices. In Appendix~\ref{sec:appA} we show how a closed form expression for a series which appear in our calculations
can be obtained by purely analytical methods.  Appendix~\ref{sec:appB} lists the notation used throughout the paper for
readers convenience.

\section{Mixed uniform measures and their joint realization}\label{sec:gce}

The generating function for the number $p(n)=\#\CP_n$ of partitions of $n$ is well known since Euler has discovered it:
\begin{equation}\label{P(q)}
P(q)=\sum_{n=0}^\infty p(n)q^n=\prod_{j=1}^\infty \frac{1}{1-q^j}.
\end{equation}
The proof (see, e.g., \cite{Andrews}) suggests the following construction first introduced by Fristedt \cite{Fristedt}.
For a number $q\in[0,1)$ introduce a measure $\mgc{q}$ on $\CP$ by 
\begin{equation}\label{Mgc}
\mgc{q}\{\la\}=P(q)^{-1} q^{N(\la)}.
\end{equation}
By definition, conditional probability measures on $\CP_n$ are uniform for any $q\in(0,1)$ and all $n\in\BN_0\eqdef \{0,1,2,\dots\}$: 
$\mgc{q}\bigl[\la\bigm|\la\in\CP_n\bigr]=M_n^{\text{uni}}\{\la\}=\frac{1}{p(n)}$.  On the other hand, due to representation \eqref{N-via-Rk}
and the product form \eqref{P(q)} of $P(q)$, part counts $C_k$, $k=1,2,\dots$, viewed as random variables on the probability space
$(\CP,\mgc{q})$ are independent and $C_k$ has the geometric distribution with success probability $1-q^k$, that is
$\mgc{q}\{\la:C_k(\la)=r\}=q^{kr}(1-q^k)$, $r\in\BN_0$.  This is similar to the construction of the
grand canonical ensemble as a mixture of microcanonical measures in statistical mechanics, as noticed
by Vershik \cite{V-FAA,V-UMN}, so following these papers we call $\mgc{q}$ a \textit{grand canonical} measure. Another 
similar probabilistic construction is known as Poissonization.
In the context of integer partitions this construction is sometimes called Fristedt's conditioning device, and it
has been used by many authors \cite{Fristedt,Pittel,Corteel,Mutafchiev,Goh-Hitczenko}. 
It has been also generalized to deal with similar measures \cite{V-FAA,Bogachev,Granovsky,Yakubovich12,BogYak}.

In a pioneering paper \cite{Fristedt} on the application of this technique to integer partitions Fristedt
proved certain properties of the grand canonical measures $\mgc{q}$ which will be instrumental below. Let for $n=1,2,\dots$
\begin{equation}\label{tn-qn}
t_n\eqdef \frac{\pi}{\sqrt{6n}},\qquad q_n\eqdef \ee^{-t_n}.
\end{equation}

\begin{lemma}[\cite{Fristedt}, Corollary 4.4 and Lemma 4.5]\label{lem:N-asymp-gc}
Governed by the probability measure $\mgc{q_n}$, 
where $q_n$ is given in \eqref{tn-qn},
the variance of $N$ is asymptotic to $v_n^2\eqdef (2\sqrt{6}/\pi)n^{3/2}$ and the difference between 
the expected value of $N$ an $n$ is $o(n^{3/4})$ as $n\to\infty$.  Moreover, for each real $x$
\begin{equation}\label{mgc-limit-thm}
\mgc{q_n}\{\la:N(\la)-n\le x v_n\}\to \frac{1}{\sqrt{2\pi}} \int_{-\infty}^x e^{-y^2/2}dy, \qquad n\to\infty,
\end{equation}
and 
\begin{equation} \label{mgc-N=n}
\mgc{q_n}\{\la:N(\la)=n\}\sim \frac{1}{\sqrt[4]{96n^3}}, \qquad n\to\infty.
\end{equation}
\end{lemma}

The idea of Fristedt's construction is to define a family of measures $\mgc{q}$ on $\CP$ so that a sequence $\mgc{q_n}$ mimics
the behavior of uniform measures on $\CP_n$.  Our approach is to enlarge the probability space and to define a 
two-parametric family of random variables $C_{q,k}$ such that  for any fixed $q\in[0,1)$ 
and for all $k=1,2,\dots$, the joint distribution of $C_{q,k}$ 
 is the same as the joint distribution of $C_k$'s when the underlying measure on $\CP$ is $\mgc{q}$.  
Since the part counts $C_k$, $k=1,2,\dots$, determine the partition, we obtain a random process on partitions 
parameterized by $q\in[0,1)$.
Moreover, our construction ensures that each $C_{q,k}$ does not decrease as a function of $q$, thus
treating $q$ as a ``time'' parameter provides a process on partitions which
is monotone (in terms of Young diagram inclusion order).  

While it is possible to do it this way, it seems more convenient to change a variable $q=\ee^{-t}$ and define processes
$\rho_k(t)$, $t\in(0,\infty)$, which are independent, do not increase in $t$ and have, for each $t>0$, the geometric distribution with success probability $1-\ee^{-kt}$.
To this end, let a double sequence $(\eps_{k,j})_{k,j\ge1}$\label{eps} of independent standard exponential random variables be defined on 
some probability space $(\Omega,\mathcal{F},\BP)$. 
Define for any $t>0$ 
\begin{equation}\label{Rk-def}
\rho_k(t)\eqdef \inf\{r\ge 0:\eps_{k,r+1}\le kt\}.
\end{equation}
We may also extend this definition to $t=\infty$ by letting $\rho_k(\infty)=0$ which makes each $\rho_k$ (left-)continuous at $\infty$ almost surely (a.s.).

\begin{proposition}\label{prop:Rks}
For a double sequence $(\eps_{k,j})_{k,j\ge1}$ of independent standard exponential random variables, let a sequence of processes
$\rho_k(t)$, $k=1,2,\dots$, $t\in(0,\infty)$, be defined by \eqref{Rk-def}. Then, for each $k$, $\rho_k(t)$ does not increase in $t$, 
for any fixed $t$ has the geometric distribution with success probability $1-\ee^{-kt}$, 
and all $\rho_k$'s are right-continuous and independent (not even for a fixed $t$ but as stochastic processes).
\end{proposition}

\begin{proof}
All claims follow immediately from the definitions.
\end{proof}

\begin{remark*}
One can avoid using double sequences by letting $\rho_k(t)=r$ for $t\in[\eps_k/(kr+k),\linebreak[3]\eps_k/(kr))$, $r=0,1,2,\dots$, $k=1,2,\dots$, 
for independent standard exponential $\eps_k$, where $\eps_k/(kr)=\infty$ for $r=0$.  This leads to the same distribution of $\rho_k(t)$,
for each fixed $t$.  However,
as it often happens, introducing additional stochasticity simplifies some arguments.
\end{remark*}

For each $t>0$ let us define a random partition $\la(t)$\label{la(t)} such that $C_k(\lambda(t))=\rho_k(t)$ for $k=1,2,\dots$.  Since all $\rho_k$'s do not increase, it follows
that $Y_{\la(s)}\subseteq Y_{\la(t)}$ for $s>t>0$.  Proposition~\ref{prop:Rks} implies that the distribution of $\la(t)$, for each fixed $t>0$,
is $\mgc{\ee^{-t}}$.  It is not defined for $t=0$.  However, for any $\tau>0$, starting from a random partition $\lambda(\tau)$ with
distribution $\mgc{\ee^{-\tau}}$ the behavior of the partition-valued process $\lambda(t)$ for $t\ge\tau$ is easily understood.

\begin{proposition}\label{prop:Markov}
For any $\tau>0$, the partition-valued random process $\la(t)_{t\ge \tau}$ is a continuous time Markov chain with initial distribution
$\mgc{\ee^{-\tau}}$ on $\CP$ defined in \eqref{Mgc}, and the generator matrix $W=(w_{\la,\mu})_{\la,\mu\in\CP}$, where
\begin{equation}
w_{\la,\mu}=\begin{cases}
-N(\la),\quad & \la=\mu,\\
\kappa,\quad & 
C_{\kappa}(\mu)<C_\kappa(\la)\text{ and }C_k(\mu)=C_k(\la)\text{ for }k\ne\kappa,\\
0,&\text{otherwise}.
\end{cases}
\end{equation}
\end{proposition}

\begin{proof}
It was already mentioned that the distribution of $\la(\tau)$ is $\mgc{\ee^{-\tau}}$.  Given
$\la(\tau)=\la^0$, the process $\la(t)_{t\ge \tau}$ is a right-continuous process on  
the \textit{finite} state space $\CP({\subseteq}\la^0)\eqdef \{\mu\in\CP:Y_\mu \subseteq Y_{\la^0} \}$.  
Given $\la(t)=\la$ for some $t\ge \tau$,
that is $\rho_k(t)=r_k$ for all $k=1,2,\dots$ with $r_k=C_k(\la)$,  the distribution of $\la(s)$, $s\ge t$, is determined
by $\eps_{k,j}$ such that $1\le j\le r_k$, and for all such $k$ and $j$ we have $\eps_{k,j}/k\ge t$. 
Then, by the memoryless property of the exponential distribution, for any $\lambda\in\CP(\subseteq\lambda^0)$ and $h>0$, 
\[
\BP[\la(t+h)=\la|\la(t)=\la]=\BP[\eps_{k,j}/k\ge h,\,k=1,2,\dots,\,1\le j\le r_k]
=\!\!\prod_{k,j:1\le j\le r_k} \!\!\!\ee^{-kh}=\ee^{-N(\la)h}
\]
and for any $\mu\in\CP(\subseteq\lambda_0)$ such that $C_{\kappa}(\mu)=r<r_\kappa$ and $C_k(\mu)=r_k$ for $k\ne\kappa$,
\begin{align*}
\BP[\la(t+h)=\mu|\la(t)=\la]
&\sim\BP[\eps_{\kappa,r+1}/\kappa\in[0,h),\,\eps_{k,j}/k\ge h\text{ for }1\le j\le r_k,(k,j)\ne(\kappa,r+1)]\\
&\sim (1-\ee^{-\kappa h})\ee^{-(N(\la)-\kappa)h}\sim \kappa h,\qquad h\downarrow0.
\end{align*}
For any other $\mu$ we have $\BP[\la(t+h)=\mu|\la(t)=\la]=O(h^2)$ as $h\downarrow 0$.
Moreover, $\la(t+h)$ is conditionally independent of $(\la(s),\tau\le s\le t)$ given $\lambda(t)$ because it is determined 
by other $\eps_{k,j}$.  Hence, by Theorem 2.8.2 of \cite{Norris},  given
$\la(\tau)=\la^0$, the process $\la(t)_{t\ge\tau}$ is a Markov chain on $\CP({\subseteq}\la^0)$ with the initial
distribution $\delta_{\la^0}$ and generator matrix $(w_{\la,\mu})_{\la,\mu\in\CP({\subseteq}\la^0)}$.  The claim follows from the fact
that $w_{\la,\mu}$ do not depend on $\la^0$.
\end{proof}

\section{Forward and backward jump chains}\label{sec:jumps}

It is possible to work with the continuous time Markov chain $\la(t)$ described in Proposition~\ref{prop:Markov}.
However it is more convenient to work with its jump chain with discrete time. For $\tau>0$ let
$\Lambda^\tau_0,\Lambda^\tau_1,\dots,\Lambda^\tau_{J(\tau)}$ be the jump chain of $\la(t)_{t\ge \tau}$.  Here $J(\tau)$ is
the (random) number of steps needed for the chain to reach its absorbing state $\Lambda^\tau_{J(\tau)}=\la(\infty)=\varnothing$.
The jump chain can be described explicitly as follows.

\begin{proposition}\label{prop:jump-forward}
The jump chain $ \Lambda^\tau_0,\Lambda^\tau_1,\dots,\Lambda^\tau_{J(\tau)}$ is the Markov chain on $\CP$ with
initial distribution $\mgc{\ee^{-\tau}}$, absorbing state $\varnothing$  and the following update rule.   Given $\Lambda^\tau_j=\la\ne\varnothing$ choose 
one of $N(\la)$ boxes in the Young diagram $Y_\la$ uniformly at random and remove a row (of some random length $\kappa$) containing the chosen box
and all rows of the same length $\kappa$ below it.  If the new diagram is $Y_\mu$, then $\Lambda^\tau_{j+1}=\mu$.   Hence, the
transition probability $p_{\la,\mu}=\kappa/n$ for all partitions $\mu$ which can be obtained from a partition $\la\in\CP_n$
by removing some number of parts $\kappa$, and zero otherwise.
\end{proposition}

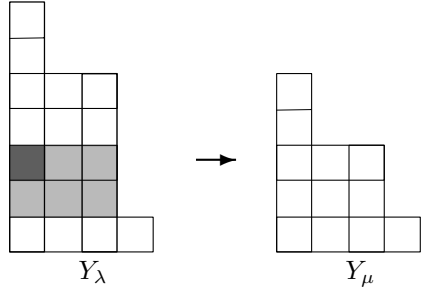
\begin{figure}
\begin{center}\begin{picture}(146,90)
\put(0,0){\includegraphics[width=0.136\textwidth]{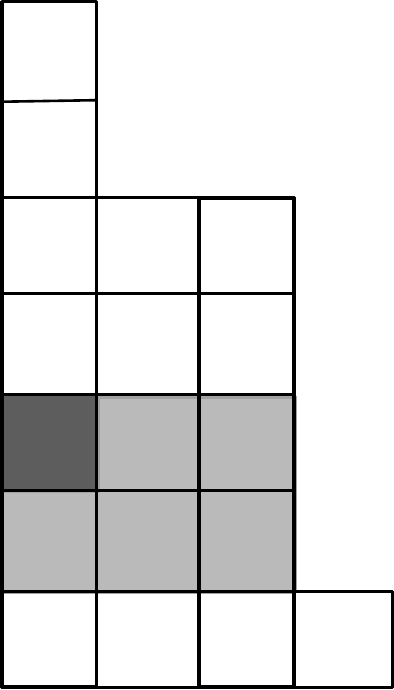}}
\put(26,-10){$Y_\la$}
\thicklines\put(70,35){\vector(1,0){15}}
\put(100,0){\includegraphics[width=0.136\textwidth]{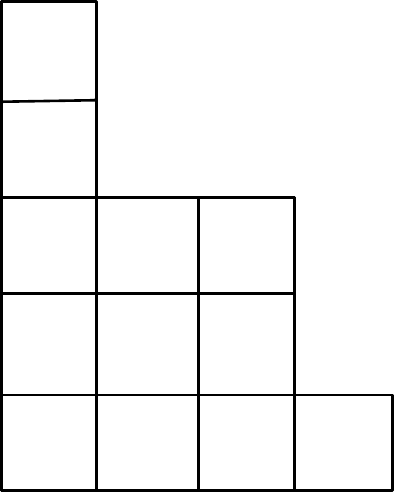}}
\put(126,-10){$Y_\mu$}
\end{picture}\end{center}
\caption{A box (dark grayed) is chosen uniformly at random in the Young diagram $Y_\la$ of $\la=(4,3,3,3,3,1,1)$
and a row containing this box is removed along with all rows of the same length below it (grayed) to obtain the Young diagram $Y_\mu$, $\mu=(4,3,3,1,1)$, in this example.}\label{fig:deletion}
\end{figure}

One step of this stochastic procedure is illustrated in Fig.~\ref{fig:deletion}.

\begin{proof} The initial distribution of $\Lambda^\tau_0$ is that of $\la(\tau)$ which is $\mgc{\ee^{-\tau}}$
according to Proposition~\ref{prop:Markov}.  By the same Proposition, the only possible jumps from the state
$\la(t)=\la$ are to states $\mu$ that are obtained from $\la$ by removing some number of  parts of the same size,
and any number of parts equal to $\kappa$ is removed with the same probability $\kappa/N(\la)$.  The description 
given in the claim clearly produces the same stochastic dynamics.
\end{proof}

We are interested in a time-reversal of the jump chain $\Lambda^\tau$.  Define a stochastic process
\begin{equation}\label{hat-Lambda}
\hLa^\tau_j = \Lambda^\tau_{J(\tau)-j},\qquad j=0,1,\dots,J(\tau).
\end{equation} 
The process $\hLa^\tau$ starts at $\hLa^\tau_0=\varnothing$ and jumps on each step from a partition $\mu$ to
a partition $\la$ which can be obtained from $\mu$ by adding some number of parts of the same size.
A key observation is that according to a result by Hunt \cite{Hunt}, $(\hLa^\tau_j)_{j=0,1,\dots,J(\tau)}$ 
is also a Markov chain with stationary transition probabilities, say  $\hp^{\,\tau}_{\mu,\la}$, which can be found from the equality
\begin{equation}\label{edge-probability}
g^\tau(\la)p_{\la,\mu}=g^\tau(\mu)\hp^{\,\tau}_{\mu,\la}
\end{equation}
where $g^\tau(\la)$ is the expected number of visits to state $\la$ for both forward and backward chains:
\begin{equation}\label{g-la}
g^\tau(\la)=\BE \sum_{j=0}^{J(\tau)}\boldsymbol{1}(\Lambda_j=\la)=\BE \sum_{j=0}^{J(\tau)}\boldsymbol{1}(\hLa_j=\la).
\end{equation}
Equation \eqref{edge-probability} means that the expected number of passes from $\la$ to $\mu$ in the forward chain
equals the expected number of passes from $\mu$ to $\la$ in the backward chain.

In order to find $g^\tau(\la)$, $\la\in\CP$, we can analyze the continuous time Markov chain $\la(t)_{t\ge \tau}$.

\begin{lemma}
For any partition $\lambda\vdash n\ge 1$, we have 
\begin{equation}\label{g-tau-la}
g^\tau(\la)=n\int_{\tau}^{\infty} \ee^{-nt}\prod_{j=1}^\infty(1-\ee^{-jt})\,\dd t.
\end{equation}
\end{lemma}

\begin{proof}
Both forward and backward chains visit each 
state no more than once.  Hence 
\[
g^\tau(\la)=\BP[\,\exists\, m\in\{0,1,\dots,J(\tau)\}:\Lambda^\tau_n=\la]=\BP[\,\exists\, t\in[\tau,\infty):\la(t)=\la].
\]
For a partition $\mu$ define $T_\mu=\sup\{t>0:\la(t)=\mu\}$, if $\la(t)\ne\mu$ for all $t>0$, then $T_\mu=-\infty$.
In this notation $g^\tau(\la)=\BP[T_\la\ge \tau]$.  For any $t>0$, 
\begin{align*}
\BP[T_\la\in[t,t+h)]&=\BP[\la(t)=\la,\la(t+h)\ne \la]+O(h^2)\\
&\sim N(\la)h\BP[\la(t)=\la]
=N(\la)h\ee^{-N(\la)t}\smash[t]{\prod_{j=1}^\infty}(1-\ee^{-jt}),\quad h\downarrow 0.
\end{align*}
The term $O(h^2)$ in the first line reflects the possibility of more than one jump in the interval $[t,t+h)$ with one of jumps into $\lambda$,
and the second line follows
because the total jump rate at state $\la$ is $N(\la)$.  Hence the continuous part of the distribution of $T_\la$ has
density $N(\la)\ee^{-N(\la)t}\prod_{j=1}^\infty(1-\ee^{-jt})$ and \eqref{g-tau-la} follows.
\end{proof} 

\begin{corollary}\label{cor:g-tau}
The probability $g^\tau(\la)$ for the chain $\la(t)_{t\ge\tau}$ to visit a partition $\la\vdash n$ can be expressed as
\begin{equation}\label{g-tau-la-as-sum}
g^\tau(\la)=\sum_{k=-\infty}^\infty \frac{(-1)^k n\ee^{-\bigl(n+\frac{k(3k+1)}2\bigr)\tau} }{n+\frac{k(3k+1)}{2}}. 
\end{equation}
One has
\begin{equation}\label{g-zero-n}
g(\la)\eqdef  \lim_{\tau\downarrow 0}g^\tau(\la)=\sum_{k=-\infty}^\infty \frac{(-1)^kn}{n+\frac{k(3k+1)}{2}}
=\frac{8\sqrt{3}\pi n\,  \sinh \left(\frac{\pi  \sqrt{24n-1}}{6}\right)}{\sqrt{24n-1} \left(2  \cosh \left(\frac{\pi  \sqrt{24n-1}}{3}\right)-1\right)}.
\end{equation}
\end{corollary}

\begin{proof}
Equation \eqref{g-tau-la-as-sum} follows from \eqref{g-tau-la} and the famous Euler pentagonal theorem \cite[Cor.~1.7]{Andrews}:
\begin{equation}\label{euler-pentagonal}
\prod_{j=1}^{\infty}(1-q^j)=\sum_{k=-\infty}^\infty (-1)^k q^{k(3k+1)/2},\qquad|q|<1.
\end{equation}
Interchanging summation and integration is justified by the dominated convergence, because 
the series  in \eqref{euler-pentagonal} converges absolutely for $|q|<1$, and taking 
$q=\ee^{-\tau}$ provides a bound. 
Changing the order between summation and limit in~\eqref{g-zero-n} is also possible by the dominated convergence, because 
taking $\tau=0$ in \eqref{g-tau-la-as-sum} gives an absolutely convergent series.  Finally, the last equality in \eqref{g-zero-n} is
proven in the Appendix.
\end{proof}

The probability for chains $\Lambda^\tau_j$ and $\hLa^\tau_j$ to visit any partition $\lambda$ depends only on its
weight $N(\la)$.  Slightly abusing the notation we write
\begin{equation}\label{gn}
g^\tau(n)\eqdef g^\tau(\la),\quad g(n)\eqdef\lim_{\tau\downarrow0}g^\tau(\lambda),\qquad \la\vdash n.
\end{equation}
The expressions for these two quantities are given in \eqref{g-tau-la-as-sum} and \eqref{g-zero-n}.

\begin{corollary}\label{cor:level-visit-prob}
As $n\to\infty$, the probability that $\la(t)$ visits $\CP_n$ is
\begin{equation}\label{visit-prob}
\gamma(n)\eqdef \BP[\,\exists\, t>0:\la(t)\in\CP_n]\sim \frac{\pi}{2\sqrt{6n}}.
\end{equation}
\end{corollary}

\begin{remark*}
This probability is much greater than the maximum over $t>0$ of the probability that $\la(t)\in\CP_n$ for 
a given $t$, which is of order $n^{-3/4}$ according to Lemma~\ref{lem:N-asymp-gc}.
\end{remark*}

\begin{proof}
By \eqref{g-tau-la} and \eqref{g-zero-n}, 
and since $g(\la)$ only depends on $n$, we get $\gamma(n)=g(n)p(n)$, where $p(n)$ is the number of partitions of $n$, as above.
From \eqref{g-zero-n} it follows that 
\begin{equation}\label{gn-asymp}
g(n)=\frac{4\sqrt{3}\pi n\,\ee^{-\pi\sqrt{24n-1}/6}}{\sqrt{24n-1}}+O\bigl(\sqrt{n}\ee^{-5\pi\sqrt{24n-1}/6}\bigr)
\sim \pi\sqrt{2n}\,\ee^{-\pi\sqrt{\frac{2n}{3}}},\quad n\to\infty.
\end{equation}
The claim follows from the famous Hardy--Ramanujan asymptotic evaluation \cite{Hardy-Ramanujan}, see also~\cite{DeSalvo}, 
\begin{equation}\label{Hardy-Ramanujan}
p(n)=\frac{2\sqrt{3}\ee^{\pi\sqrt{24n-1}/6}}{24n-1}\Bigl(1-\frac{6}{\pi\sqrt{24n-1}}\Bigr)+O\Bigl(\frac{\ee^{\pi\sqrt{n/6}}}{n}\Bigr)
\sim\frac{1}{4\sqrt{3}\,n}\ee^{\pi\sqrt{\frac{2n}{3}}},\quad n\to\infty. \qedhere
\end{equation}
\end{proof}

This allows us to describe a limiting process, as $\tau\downarrow0$, for the backward chain $\hLa^\tau$, as a 
Markov chain with explicit transition probabilities.  For notational convenience, let us
define $\hLa^\tau_j:=\hLa^\tau_{J(\tau)}$ for $j>J(\tau)$. (This is not a Markov process, because the transition 
probabilities change at random time $J(\tau)$ to $\hp^{\,*}_{\la,\la}=1$, but up to time $J(\tau)$ the behavior is Markovian
with transition probabilities given by $\hp^{\,\tau}_{\mu,\la}$ defined in \eqref{edge-probability}.)

\begin{thm}\label{th:hat-Lambda}
For each fixed $m$, random vectors $(\hLa^\tau_0,\hLa^\tau_1,\dots,\hLa^\tau_m)$
converge almost surely, as $\tau\downarrow0$, to a random vector $(\hLa_0,\hLa_1,\dots,\hLa_m)$.
The latter is the jump process of the process $\la(t)$ viewed backward from time $t=\infty$ and 
is a Markov chain which starts at $\hLa_0=\varnothing$ and has stationary
transition probabilities
\begin{equation}\label{hatp}
\hp_{\mu,\la}\eqdef\begin{cases}
\dfrac{\kappa g({n+\kappa r})}{(n+\kappa r)g(n)},  \quad &
\begin{aligned}&C_\kappa(\la)=C_\kappa(\mu)+r\text{ for some }\kappa,r\ge1,\\[-0.6ex]
&C_k(\la)=C_k(\mu)\text{ for }k\ne\kappa,\end{aligned}\\[0.6ex]
0, & \text{otherwise},
\end{cases}
\end{equation}
where $n=N(\mu)$ is the weight of the partition $\mu$ and $g(n)$ is defined in \eqref{gn},~\eqref{g-zero-n}.
\end{thm}

\begin{proof}
Let $(\hLa_j)_{j\ge 0}$ be the jump process of the chain $(\la(t))_{t>0}$ viewed backward from time $t=\infty$, 
that is $\hLa_0=\la(\infty)=\varnothing$ and further values are defined recurrently: if $T_j=\inf\{t:\la(t)=\hLa_j\}$ for 
some $j=0,1,2,\dots$ the next value is $\hLa_{j+1}=\la(T_j-0)$.
By construction, given that the number of jumps $J(\tau)$ of the chain $(\la(t))_{t\ge\tau}$ is 
at least $m$, we have $\hLa^\tau_j=\hLa_j$ for $j=0,1,\dots,m$.  Hence the almost sure convergence follows from
\begin{equation}\label{M-tau-ge-M}
\BP[J(\tau)\ge m]\to 1,\qquad\tau\downarrow0,
\end{equation}
for any fixed $m$, and \eqref{M-tau-ge-M} holds true because $J(\tau)$ is
not less than the number of
lower records greater than $\tau$ in the infinite sequence
$\eps_{1,1},\eps_{1,2},\dots$.

To find the transition probabilities of $(\hLa_j)$ note that
on the event $\{J(\tau)\ge m\}$ the process $(\hLa^\tau_j)_{j=0,\dots,m}$ is the
Markov chain.  Its transition probabilities are $\hp^{\,\tau}_{\mu,\la}/\BP[J(\tau)\ge m]$, where $\hp^{\,\tau}_{\mu,\la}$ is defined in \eqref{edge-probability}. 
As $\tau\downarrow0$, they converge to $\hp_{\mu,\la}$ given by \eqref{hatp}, 
according to Corollary~\ref{cor:g-tau} and \eqref{M-tau-ge-M}.
\end{proof}

The chain $(\hLa)_{j\ge0}$ evolves by inserting a rectangle of some random size $K_j\times R_j$ into the Young diagram $Y_{\hLa_j}$.
The whole sequence $(\hLa_j)_{j\ge 0}$ can be reconstructed from a sequence of pairs $(K_j,R_j)_{j\ge0}$.   However,
this sequence of pairs is not Markovian, because the distribution of $(K_j,R_j)$ depends on the partition $\hLa_j$. 
Theorem~\ref{th:hat-Lambda} shows that actually it depends just on the weight $N(\hLa_j)$ of the partition,
so the sequence of triples $(N(\hLa_j),K_j,R_j)_{j\ge 0}$ is the Markov chain. 
Another Markov chain embedded in $\hLa$ is the sequence $(N(\hLa_j))_{j\ge 0}$, but it does not determine $\hLa$.
We proceed with investigating these two chains.

\begin{proposition}\label{prop:hN-increasing}
The process $(N(\hLa_j))_{j\ge 0}$ is an increasing Markov chain on non-negative integers which starts 
at $N(\hLa_0)=0$ and has stationary transition probabilities 
\begin{equation}\label{hq}
\hq_{n,n+h}\eqdef\frac{g(n+h)\sigma_1(h)}{(n+h)g(n)},\qquad h=1,2,\dots,
\end{equation}
where $\sigma_1(h)=\sum_{r|h}r$ is the sum of divisors function \textup{(}$r|h$ means ``$r$ divides $h$''\textup{)}
 and $g(n)$ is defined in \eqref{gn},~\eqref{g-zero-n}.
\end{proposition}

\begin{proof}
If $N(\hLa_j)=n$ for some $j,n\ge 0$, the next partition $\hLa_{j+1}$ can have any weight $n+h$, $h\ge 1$,
if one inserts a rectangle $k\times r$ into $Y_{\hLa_j}$, for some $k,r$ such that $kr=h$.  Hence, from \eqref{hatp} we
obtain, for $h\ge 1$,
\[
\BP[N(\hLa_{j+1})=n+h|N(\hLa_j)=n]=\sum_{r|h}\frac{(h/r)g(n+h)}{(n+h)g(n)} 
=\frac{g(n+h)}{(n+h)g(n)}\sigma_1(h),
\]
because when $r$ runs through the set of divisors of $h$, $h/r$ runs through the same set.
\end{proof}

\begin{corollary}\label{cor:qnm}
For any $\tau>0$, the process $(N(\Lambda_j^\tau))_{j=0,\dots,J(\tau)}$ is a Markov chain  
with the initial distribution $\BP[N(\Lambda_0^\tau)=n]=p(n)\ee^{-n\tau}P(\ee^{-\tau})^{-1}$ and stationary 
transition probabilities
\begin{equation}\label{qnm}
q_{n,m}:=\BP[N(\Lambda_{j+1}^\tau)=m|N(\Lambda_{j}^\tau)=n]=\frac{p(m)\sigma_1(n-m)}{np(n)},\qquad 0\le m<n,
\end{equation}
where $\sigma_1(h)$ is the sum of divisors function and $p(n)=|\CP_n|$ is the number of partitions of~$n$.
\end{corollary}

\begin{proof}
The process $(N(\Lambda_j^\tau))_{j=0,\dots,J(\tau)}$ is a time-inverse of the process
$(N(\hLa_j))_{j=0,\dots,J(\tau)}$: $N(\Lambda_j^\tau)=N(\hLa_{J(\tau)-j})$.  (Note that the absorption
time $J(\tau)$ depends on the trajectory $(\hLa_j)_{j\ge 0}$, however it does not affect the next argument.)
By the result of Hunt \cite{Hunt} this time-inverse is a Markov chain, and the transition probabilities satisfy
\begin{equation}\label{q-hq}
\gamma(n)q_{n,m}=\gamma(m)\hq_{m,n},\qquad 0\le m<n,
\end{equation}
where $\gamma(n)=p(n)g(n)$ (see \eqref{visit-prob}), which gives \eqref{qnm}.
\end{proof}

The fact that $(q_{n,m})_{m=0,\dots,n-1}$ is a probability vector for each $n\in\BN$ is equivalent to a well-known recursion
\begin{equation}\label{npn-recursion}
np(n)=\sum_{m=0}^{n-1}p(m)\sigma_1(n-m),
\end{equation}
see, e.g., \cite[Ex.~3.77]{Aigner}.  We were unable to find the claim that  $(\hq_{n,m})_{m>n}$ is a probability
distribution, that is 
\begin{equation}\label{gn-recurstion}
\sum_{m=n+1}^\infty \frac{g(m)\sigma_1(m-n)}{m}=g(n)
\end{equation} 
for each~$n\in\BN_0$, in the literature.  Let us rewrite \eqref{gn-recurstion} in more explicit form. To this end
introduce a notation
\begin{equation}\label{Ankr}
A_n(h)=\tfrac{\pi}{6}\sqrt{24n+24h-1}.
\end{equation} 
From \eqref{g-zero-n}, using this notation we have
\begin{equation}\label{gnkr}
g(n+h)=\frac{4 \sqrt{3}\pi^2 (n+h)  \sinh \left(A_n(h)\right)}{3A_n(h) \left(2  \cosh \left(2A_n(h)\right)-1\right)}.
\end{equation}
Hence an infinite sequence of relations between sum of divisors function values holds: for each $n=0,1,2,\dots$,
\begin{equation}\label{divisor-function}
\sum_{h=1}^\infty \frac{  \sinh \left(A_n(h)\right)\,\sigma_1(h)}{A_n(h) \left(2  \cosh \left(2A_n(h)\right)-1\right)}
=\frac{ n\,  \sinh \left(A_n(0)\right)}{A_n(0) \left(2  \cosh \left(2A_n(0)\right)-1\right)}.
\end{equation}
(For $n=0$ the right-hand side is $3g(0)/(4\sqrt3\pi^2)=\sqrt{3}/(4\pi^2)$ which is also the limit of the right-hand side as $n\to 0$.)   
One can think of equations \eqref{divisor-function}
as of an infinite linear system for variables $(\sigma_1(h))_{h\ge 1}$.  The matrix of this system is a Hankel matrix. 
We could not manage to diagonalize this matrix and obtain an explicit expression for $\sigma_1(h)$.

\smallskip

If $N(\hLa_j)=n$, the distribution of $K_j$ and $R_j$ can be read from \eqref{hatp}.  It seems that 
simple explicit expressions for the marginal distributions do not exist, however their asymptotics can be calculated.  Here and below we
use notation $\boldone_A$ for the indicator of an event $A$.

\begin{proposition}\label{prop:kr-limits}
As $n$ grows unboundedly, we have the following asymptotic behavior of $K_j$ and $R_j$ given $N(\hLa_j)=n$:
\begin{enumerate}
\item[\textup{(i)}]  $K_j$ grows proportionally to $\sqrt{n}$: for any $x>0$
\begin{equation}\label{k-limit}
\lim_{n\to\infty}\BP\bigl[K_j\le x\sqrt{n}\big|N(\hLa_j)=n\bigr]=\int_0^x \frac {y\, \ee^{-\pi y/\sqrt{6}}dy}{1-\ee^{-\pi y/\sqrt{6}}}\,;
\end{equation}
\item[\textup{(ii)}] for any $x>0$, given $K_j\sim x\sqrt{n}$, 
$R_j$ has asymptotically geometric distribution on $\{1,2,\dots\}$ 
with parameter $1-\ee^{-\pi x/\sqrt{6}}$\textup{: }for any sequence $k_n$ such that $\smash[b]{\lim\limits_{n\to\infty}} k_n/\sqrt{n}=x$ and for each $r=1,2,\dots,$
\begin{equation}\label{r-limit}
\lim_{n\to\infty}\BP\bigl[R_j=r\big|K_j=k_n,N(\hLa_j)=n\bigr]= \bigl(1-\ee^{-\pi x/\sqrt{6}}\bigr)\ee^{-\pi x(r-1)/\sqrt{6}};
\end{equation}
\item[\textup{(iii)}] in one step, the mean Young diagram area increases 
proportionally to $\sqrt{n}$:
\begin{equation}\label{EKR}
\lim_{n\to\infty}\frac{\BE\bigl[K_jR_j\big|N(\hLa_j)=n\bigr]}{\sqrt{n}}=\frac{2\sqrt{6}}{\pi}
\end{equation} 
in accordance with \eqref{visit-prob}, moreover, for any $\nu=1,2,\dots$,
\begin{equation}\label{EKRm}
\BE\bigl[(K_jR_j)^\nu\big|N(\hLa_j)=n\bigr]=O(n^{\nu/2}),\qquad n\to\infty;
\end{equation} 
\item[\textup{(iv)}] for any $x>0$, 
\begin{equation}\label{ER}
\lim_{n\to\infty}\BE\bigl[R_j\boldone_{\{K_j\ge x\sqrt{n}\}}\big|N(\hLa_j)=n\bigr]=h(x),
\end{equation}
where 
\begin{equation}\label{h-def}
h(x):=\frac{\sqrt{6}\,x\,\ee^{-\pi x/\sqrt{6}}}{\pi(1-\ee^{-\pi x/\sqrt{6}})}-\frac{6\log(1-\ee^{-\pi x/\sqrt{6}})}{\pi^2}\,.
\end{equation}
\end{enumerate}
\end{proposition} 

\begin{proof} 
Rewrite the conditional distribution of $(K_j,R_j)$ given 
$N(\hLa_j)=n$ in terms of $A_n(kr)$:
\begin{align}\label{kr-given-n}
\notag
\BP[K_j=k,R_j=r|N(\hLa_j)=n]&=\dfrac{k g({n+kr})}{(n+kr)g(n)}\\
&=\frac{ k A_n(0) \sinh \left(A_n(kr)\right)\left(2  \cosh \left(2A_n(0)\right)-1\right) }{nA_n(kr)   \sinh \left(A_n(0)\right)\left(2  \cosh \left(2A_n(kr)\right)-1\right)}.
\end{align}
The ratio of hyperbolic functions is
\begin{equation}\label{hyperbolic-ratio}
\frac{\sinh \left(A_n(kr)\right)\left(2  \cosh \left(2A_n(0)\right)-1\right) }{ \sinh \left(A_n(0)\right)\left(2  \cosh \left(2A_n(kr)\right)-1\right)}
=\exp(A_n(0)-A_n(kr))\bigl(1+O(\ee^{-\pi\sqrt{\frac{2n}{3}}})\bigr)
\end{equation}
as $n\to\infty$, where $O(\cdot)$ is uniform in $k,r$.
We have
\[
\frac{A_n(kr)}{A_n(0)}=\frac{\sqrt{24n-1+24kr}}{\sqrt{24n-1}}=\sqrt{1+\frac{kr}{n-\tfrac{1}{24}}}.
\]
In view of elementary inequalities $1+x/2-x^2/2\le \sqrt{1+x}\le 1+x/2$ valid for $x\ge -1$, and 
$\frac{1}{n}\le \frac{1}{n-1/24}\le \frac1n(1+\frac{1}{23n})\le \frac2n$ valid for $n\ge 1$,
we obtain bounds
\[
1+\frac{kr}{2n}-2\Bigl(\frac{kr}{n}\Bigr)^2
\le 1+\frac{kr}{2n-\frac1{12}}-\frac{1}{2}\Bigl(\frac{kr}{n-\frac1{24}}\Bigr)^2
\le\frac{A_n(kr)}{A_n(0)}\le 1+\frac{kr}{2n-\frac1{12}}\le 1+\frac{kr}{2n}\Bigl(1+\frac{1}{23n}\Bigr).
\]
Thus 
\[
-\frac{krA_n(0)}{2n}\Bigl(1+\frac{1}{23n}\Bigr)
\le A_n(0)-A_n(kr)
\le -\frac{krA_n(0)}{2n} +2A_n(0)\Bigl(\frac{kr}{n}\Bigr)^2.
\]
Similarly, we obtain a bound for the opposite ratio:
\begin{equation}\label{An0/Ankr}
1-\frac{kr}{n}\le 1-\frac{kr}{2n-\frac{1}{12}}\le \frac{A_n(0)}{A_n(kr)}\le 1.
\end{equation}
This way we obtain a two-sided bound: for $n$ large enough
\begin{multline}\label{two-sided-bounds}
\frac{k}{n}\Bigl(1-\frac{kr+1}{n}\Bigr)\exp\Bigl(-\frac{krA_n(0)}{2n}\Bigl(1+\frac{1}{23n}\Bigr)\Bigr) 
\le \dfrac{k g({n+\kappa r})}{(n+k r)g(n)}\\ 
\le \frac{k}{n}\Bigl(1+\frac{1}{n}\Bigr)\exp\Bigl(-\frac{krA_n(0)}{2n} +2A_n(0)\Bigl(\frac{kr}{n}\Bigr)^2\Bigr)
\end{multline}
where $O(\ee^{-\pi\sqrt{\frac{2n}{3}}})$  from \eqref{hyperbolic-ratio} is replaced by $\frac1n$.

Our first aim is to sum both bounds over $r=1,2,\dots$, but for $r$ large enough the upper bound in
\eqref{two-sided-bounds} becomes loose.  
Note that
\[
s\leq \sqrt{24 (n+k r)-1}<s+1\quad 
\text{ for }\frac{-24 n+s^2+1}{24 k}\leq r<\frac{-24 n+s^2+2 s+2}{24 k}.
\]
Hence combining \eqref{kr-given-n}, \eqref{hyperbolic-ratio} and the upper bound in \eqref{An0/Ankr} we obtain, for $s_0>\sqrt{24n-1}$,
\begin{equation}\label{r-tail}
\sum_{r\ge \frac{-24 n+s_0^2+1}{24 k}}\dfrac{k g({n+\kappa r})}{(n+k r)g(n)}
\le 2\frac{k}{n}\ee^{A_n(0)}\sum_{s\ge s_0}\frac{1+2s}{24k}\ee^{-\pi s/6} 
\le \frac{s_0}{n}\ee^{A_n(0)-\pi s_0/6}.
\end{equation}
Take $s_0=\sqrt{24(n+n^\alpha)-1}$ for some $\alpha>\tfrac12$.  Then 
\[
\tfrac{6}{\pi}A_n(0)- s_0=\sqrt{24n-1}-\sqrt{24(n+n^\alpha)-1}=\frac{-24n^\alpha}{\sqrt{24n-1}+\sqrt{24(n+n^\alpha)-1}}\le -n^{\alpha-1/2}
\]
and the right-hand side of \eqref{r-tail} is bounded above by $\sqrt{\frac{48}{n}}\ee^{-n^{\alpha-1/2}/2}=o(n^{-1/2})$, $n\to\infty$.
The sum over  $r<\frac{-24 n+s_0^2+1}{24 k}=n^\alpha/k$ can be restricted from both sides by expressions which have the same asymptotics as $n\to\infty$
for any $\alpha\in(\tfrac12,1)$:
\begin{multline}\label{two-sided-bounds-sum}
\frac{k}{n}\Bigl(1-\frac{2n^\alpha}{n}\Bigr)
\frac{\exp\Bigl(-\frac{kA_n(0)}{2n}\Bigl(1+\frac{1}{23n}\Bigr)\Bigr)-\exp\Bigl(-\frac{n^\alpha A_n(0)}{2n}\Bigl(1+\frac{1}{23n}\Bigr)\Bigr)}
{1-\exp\Bigl(-\frac{kA_n(0)}{2n}\Bigl(1+\frac{1}{23n}\Bigr)\Bigr)} 
\le \!\!\!\sum_{r=1}^{\frac{-24 n+s_0^2+1}{24 k}}\!\!\dfrac{k g({n+\kappa r})}{(n+k r)g(n)}\\
\le  \frac{k}{n}\Bigl(1+\frac{1}{n}\Bigr)
\frac{\exp\Bigl(-\frac{kA_n(0)}{2n}\bigl(1-2n^{\alpha-1}\bigr)\Bigr)-\exp\Bigl(-\frac{n^\alpha A_n(0)}{2n}\bigl(1-2n^{\alpha-1}\bigr)\Bigr)}
{1-\exp\Bigl(-\frac{kA_n(0)}{2n}\bigl(1-2n^{\alpha-1}\bigr)\Bigr)} .
\end{multline}
Taking some $\alpha\in(\tfrac{1}{2} ,1)$  and combining these bounds with \eqref{r-tail} we see that uniformly in~$k$
\begin{equation}\label{k-asympotics}
\BP[K_j=k|N(\hLa_j)=n]=\sum_{r=1}^\infty \dfrac{k g({n+\kappa r})}{(n+k r)g(n)}\sim \frac{k}{n}
\frac{\exp\Bigl(-\frac{kA_n(0)}{2n}\Bigr)}
{1-\exp\Bigl(-\frac{kA_n(0)}{2n}\Bigr)} ,\qquad n\to\infty.
\end{equation}

The probability that $K_j\le x\sqrt{n}$ is then the sum of expressions \eqref{k-asympotics}
(or, more exactly,  is bounded by the sums of upper and lower bounds in \eqref{two-sided-bounds-sum}) over $k\le x\sqrt{n}$.  For 
$k\sim y\sqrt{n}$ one has $kA_n(0)/(2n)\sim \pi y/\sqrt{6}$ as $n\to\infty$ uniformly in $y\le x$ for any fixed $x$.  
This gives the first claim \eqref{k-limit} due to the integral approximation of the sum.

The second claim \eqref{r-limit} follows by combining \eqref{k-asympotics} with \eqref{two-sided-bounds}.

For the third claim \eqref{EKR} we have to multiply \eqref{two-sided-bounds} by $kr$ and sum over all $k,r\ge 1$.
For the sum over $kr\ge n^\alpha$, for $\alpha\in(1/2,3/4)$ the exponential bound $O\bigl(\ee^{-n^{\alpha-1/2}/2}\bigr)$ is obtained similarly
to \eqref{r-tail}. For other values of $k$ and $r$ we first sum over $r$ 
using the asymptotic relation
\begin{equation}\label{sum-rx^r}
\sum_{r=1}^{n^\alpha/k}r \ee^{-r k a_n}=\frac{\ee^{-k a_n}}{(1-\ee^{-k a_n})^2} +O\bigl(\ee^{-n^{\alpha-1/2}/2}\bigr),\qquad n\to\infty,
\end{equation}
with $a_n=A_n(0)/(2n)\sim \pi/\sqrt{6n}$, and then use the integral approximation to calculate
\[
\BE\bigl[K_jR_j\big|N(\hLa_j)=n\bigr]=\bigl(1+O(\tfrac{1}{\sqrt{n}})\bigr)\int_0^\infty \frac{y^2\ee^{-\pi y/\sqrt{6}}dy}{(1-\ee^{-\pi y/\sqrt{6}})^2}\to \frac{2\sqrt{6}}{\pi},
\qquad n\to\infty.
\]
To show \eqref{EKRm} we should multiply \eqref{two-sided-bounds} by $(kr)^\nu$ instead. Then we can bound a sum over $kr\ge n^{\alpha}$
for $\alpha\in(1/2,3/4)$ similarly to \eqref{r-tail}, use the inequality $r^\nu\le r(r+1)\dots(r+\nu-1)$ and the closed form summation formula
\[
\sum_{r=1}^\infty r(r+1)\dots(r+\nu-1)\ee^{-kra_n} =\frac{\nu!\,\ee^{-ka_n}}{(1-\ee^{-ka_n})^{m+1}} 
\] 
to obtain, for $n$ large enough, a bound 
\[
\BE[(K_jR_j)^\nu|N(\hLa_j)=n]\le \sum_{k=1}^\infty \frac{2k^{\nu+1}\nu!\,\ee^{-ka_n}}{n(1-\ee^{-ka_n})^{\nu+1}} 
\le 4\nu!\,n^{\nu/2}\int_0^\infty \frac{y^{\nu+1}\ee^{-\pi y/\sqrt{6}}}{(1-\ee^{-\pi y/\sqrt{6}})^{\nu+1}}dy
\]
where the integral converges.

Finally, to obtain the fourth claim \eqref{ER} we multiply \eqref{two-sided-bounds} just by $r$
and sum over all $r\ge 1$ and $k\ge x\sqrt{n}$.  The tail sum over $r> n^\alpha/k$ is bounded
similarly to \eqref{r-tail}, while for the main part we 
again first use \eqref{sum-rx^r} and then the integral approximation.  This is a straightforward calculation, so the details are omitted. 
\end{proof} 

\section{Random partitions and Poisson point processes}\label{sec:Ppp}

It is well known (see, e.g.,~\cite{Stam}) that for a fixed $k$ lower records of the sequence $\eps_{k,1}/k,\eps_{k,2}/k,\dots$ form a
Poisson point process $\Pi_k$\label{Pik-def} on $\BR_+$ with the density
\begin{equation}\label{fk-def}
f_k(t) = \frac{k\ee^{-tk}}{1-\ee^{-tk}}
\end{equation}
of the intensity measure.
The process $\rho_k$ jumps precisely at points of this Poisson point process.
The jump at point $t$ can be interpreted as a mark of this point, see~\cite[Ch.~5]{Kingman}.
It has geometric distribution on $\{1,2,\dots\}$ with success probability $1-\ee^{-tk}$. 
Indeed,  $\rho_k(t)$ jumps at point $t$ from $r+j$ to $r$ if and only if there is a lower record $\eps_{k,r+1}=kt$
and the next lower record is $\eps_{k,r+j+1}$,   that is when $\eps_{k,1},\dots,\eps_{k,r},\eps_{k,r+2},\dots,\eps_{k,r+j}>kt$,
$\eps_{k,r+1}=kt$ and $\eps_{k,r+j+1}\le kt$. Hence, given a downward jump $R$
occurs at time $t$, it has size $j=1,2,\dots$ with probability 
$(1-\ee^{-tk})\,\ee^{-tk(j-1)}$.  The distribution of the jump $R$ depends on the point of jump $t$, but the value of a jump
$R=\rho_k(t-0)-\rho_k(t)$ is independent of $\rho_k(t)=r$ because it is determined by the next independent exponentials $\eps_{k,r+2},\eps_{k,r+3},\dots$.

The jumps of $\lambda(t)$ occur when either of chains $\rho_k$ jumps, that is at points of a Poisson point process
$\Pi=\cup_{k\ge 1}\Pi_k$\label{Pi-def}.  Since the point processes $\Pi_k$ are independent, $\Pi$ is also a Poisson point
process with the density 
\begin{equation}\label{f-def}
f(t) \eqdef \sum_{k=1}^\infty f_k(t)
\end{equation}
of the intensity measure.
Given $\Pi\ni t$, there is a.s.\ just one process $\Pi_k$ which jumps at $t$. The conditional probability
that it is $\Pi_\kappa$ for some specific $\kappa=1,2,\dots$ is $f_\kappa(t)/f(t)$, and $\kappa$ is 
independent of all other point of $\Pi$ (\cite[Ch.~5]{Kingman}).  Hence the random variable $K$, the size of
parts which are removed from $\lambda(t)$ at moment $t\in\Pi$, can be also treated as an independent mark of 
a point $t$ in the process $\Pi$.  Hence we can consider the Poisson point process $\Pi$ with intensity density $f(t)$
which consists of points $T_0>T_1>T_2>\dots$. For each $j=0,1,2,\dots$, the point $T_j$ is marked by a random pair $(K_j,R_j)$ of positive integers,
which is independent of $\Pi$ and has the joint distribution 
\begin{equation}\label{KjRj-joint}
\BP[K_j=k,R_j=r]=k\ee^{-tkr}/f(t),\qquad k=1,2,\dots,\quad r=1,2,\dots,
\end{equation}
given $T_j=t$, for any $j$.
 
The jump of $N(\la(t))$ at a point $T_j=t$ of $\Pi$ is then $K_jR_j$.   The distribution 
of the jump is 
\[
\BP[K_jR_j=h]=\sum_{k|h}\frac{k\ee^{-th}}{f(t)}=\frac{\sigma_1(h)\ee^{-th}}{f(t)}
\]
given $T_j=t$.  The fact that this is a probability distribution is equivalent to a well known 
Lambert series for the generating function of $\sigma_1(h)$:  returning to the variable $q=\ee^{-t}$ we obtain
\[
\sum_{h=1}^\infty \sigma_1(h)q^h=\sum_{k=1}^\infty \frac{kq^k}{1-q^k}.
\]
This shows that $f(t)$ is a Laplace generating function for $\sigma_1(h)$:
\begin{equation}\label{f-via-sigma1}
f(t)=\sum_{h=1}^{\infty}\sigma_1(h)\ee^{-th}.
\end{equation}

\smallskip 

Let us introduce tails of intensity measures
\begin{align}\label{Fk(t)}
F_k(t)&=\int_{t}^\infty \!f_k(s)\,ds=-\log(1-\ee^{-kt}),\\ 
\label{F(t)}
F(t)&=\int_{t}^{\infty} \!f(s)\,ds=-\sum_{k=1}^\infty \log(1-\ee^{-kt})=\sum_{h=1}^\infty \frac{\sigma_1(h)}{h}\ee^{-th}.
\end{align}
The representation of the partition growth process in terms of marked Poisson point processes has at least two consequences.  The
first concerns the distribution of the number of jumps $J(\tau)$ of the chain $\Lambda^\tau_j$ and follows immediately from 
the description of jump points as a Poisson point process.

\begin{proposition}
For any $\tau>0$, the number $J(\tau)$ of jumps that the Markov chain $(\la(t))_{t\ge \tau}$ needs to reach its absorbing state $\lambda(\infty)=\varnothing$ 
has the Poisson distribution with parameter $F(\tau)$.\qed
\end{proposition}

The second consequence is a possibility to introduce another time change so that the random partition grows
as time increases, and growth happens at jump times of the standard Poisson process, thus providing intuitively more clear picture.  
To this end let $s=F(t)$\label{s} and 
define $S_j=F(T_j)$ for $j=0,1,2,\dots$.  Note that points $T_j$ are ordered by decrease, so $S_0<S_1<\dots$ a.s., because $F$ 
decreases from infinity at $0$ to zero at $\infty$.
It is well known (see, e.g., \cite[Ch.~4]{Kingman}) and easy to check that $S_j$ are points of the standard Poisson point process $\tPi$\label{tPi-def} (with unit intensity).
The time $t$ of the original process $\Pi$ is connected to the time $s$ of the transformed process $\tPi$ via $t=F^{-1}(s)$.
Hence, rewriting \eqref{KjRj-joint} in the new time, we see that the distribution of a mark $(\tK_j,\tR_j)$ of a point $S_j=s$ is 
\begin{equation}\label{tilde-KR}
\BP[\tK_j=k,\tR_j=r]=k\ee^{-F^{-1}(s)kr}/f(F^{-1}(s)),\qquad k=1,2,\dots,\quad r=1,2,\dots.
\end{equation}
We can reconstruct the random partition $\tla(s)=\lambda(F^{-1}(s))$\label{tla(s)} as a partition with $\trho_k(s)$ parts $k$, where
$\trho_k$ is determined from the standard Poisson point process $\tPi=\{S_0,S_1,\dots\}$ with points $S_j$ marked by pairs $(\tK_j,\tR_j)$ as
\begin{equation}\label{trhok-def}
\trho_k(s)=\sum_{j:S_j\le s}\tR_j\boldone_{\{\tK_j=k\}}.
\end{equation}

It is easy to find the asymptotics of $f(t)$ and $F(t)$ as $t\downarrow 0$, say, by using integral approximation
and monotonicity:
\begin{equation}\label{f-F-t-asymptotics}
f(t)=\frac{\pi^2}{6t^2}+O(\tfrac1t),\qquad F(t)=\frac{\pi^2}{6t}+O(\log\tfrac1t),\qquad t\downarrow 0.
\end{equation}
From this asymptotic relations one readily finds that
\begin{equation}\label{f-F-s-asymptotics}
F^{-1}(s)=\frac{\pi^2}{6s}+O(\tfrac1{s^2}\log s), \qquad f(F^{-1}(s))=\frac{6s^2}{\pi^2}+O(s\log s),\qquad s\to\infty.
\end{equation}
This gives an asymptotics for the distribution of $(\tK_j,\tR_j)$ given the marked point is $S_j=s$ as $s\to\infty$: uniformly in $k,r\in\BN$, 
\begin{align}\label{tKtR-joint}
\BP[\tK_j=k,\tR_j=r] &= \frac{1}{s}\,\frac{\pi^2}{6}\,\frac{k}{s}\exp\Bigl(-\frac{\pi^2}{6}\,\frac{kr}{s}\Bigr)\bigl(1+O(\tfrac{\log s}{s})\bigr),\\
\BP[\tK_j=k] &= \frac{1}{s}\,\frac{\pi^2}{6}\,\frac{k}{s}\frac{\exp\bigl(-\frac{\pi^2}{6}\,\frac{k}{s}\bigr)}{1-\exp\bigl(-\frac{\pi^2}{6}\,\frac{k}{s}\bigr)}\bigl(1+O(\tfrac{\log s}{s})\bigr),\\
\BP[\tR_j=r|\tK_j=k]&=\Bigl(1-\exp\Bigl(-\frac{\pi^2}{6}\,\frac{k}{s}\Bigr)\Bigr)\exp\Bigl(-\frac{\pi^2}{6}\,\frac{k(r-1)}{s}\Bigr)\bigl(1+O(\tfrac{\log s}{s})\bigr).
\end{align}
This way, by a straightforward asymptotic evaluation, we obtain the following counterpart
of Proposition \ref{prop:kr-limits}(i,\,ii).

\begin{proposition}\label{prop:tktr-limits}
The distribution of a (scaled) mark $(\tK_j/s,\tR_j)$ of a point $s=S_j\in\tPi$ converges, as $s\to\infty$, to
the distribution of a pair $(K,R)$ such that $K$ is absolutely continuous with respect to the Lebesgue measure,
\begin{equation}\label{tk-limit}
\BP[K<x]=\frac{\pi^2}{6}\int_0^x\frac{y\,\ee^{-\pi^2y/6}}{1-\ee^{-\pi^2y/6}}\, dy,\qquad x>0,
\end{equation}
and given $K$, $R$ has the geometric distribution
\begin{equation}\label{tr-limit}
\BP[R=r|K]=\Bigl(1-\exp\Bigl(-\frac{\pi^2K}{6}\Bigr)\Bigr)\exp\Bigl(-\frac{\pi^2(r-1)K}{6}\Bigr),\quad r=1,2,\dots.
\end{equation}
\end{proposition}

Note that for any $t>0$
\[
f(t)=\BE N(\la(t)).
\]
So, according to \eqref{f-F-t-asymptotics} and \eqref{f-F-s-asymptotics}, one should take
\[
t_n\sim \frac{\pi}{\sqrt{6n}} \quad\text{and}\quad s_n\sim\pi \sqrt{\frac{n}{6}},\qquad n\to\infty, 
\]
in order to have 
\[
\BE N(\la(t_n))=\BE N(\tla(s_n))=n.
\]
This is asymptotically the choice made in Lemma~\ref{lem:N-asymp-gc}. 
Let us take $t_n=\pi/\sqrt{6n}$ as in \eqref{tn-qn}, and 
\begin{equation}\label{sn}
s_n\eqdef F(t_n)=\pi\sqrt{\frac{n}{6}}+O(\log n),\qquad n\to\infty.
\end{equation}
Then, according to Lemma~1, we have $|N(\tla(s_n))-n|=O(n^{3/4})$ with 
probability approaching 1 as $n\to\infty$.   Taking $s=s_n$ in Proposition~\ref{prop:tktr-limits}
we obtain a scaling different from that in Proposition \ref{prop:kr-limits} (by a factor $\sqrt{6}/\pi$) 
which explains the difference in formulas \eqref{k-limit}--\eqref{r-limit} and \eqref{tk-limit}--\eqref{tr-limit}.

\begin{remark*}
While the proof of Proposition \ref{prop:tktr-limits} is computationally much simpler than
the proof of Proposition \ref{prop:kr-limits} and they give consistent results, it does not
seem easy to derive Proposition \ref{prop:kr-limits} from Proposition \ref{prop:tktr-limits} directly.
\end{remark*} 

\section{The limit shape}\label{sec:ls}

Recall that we have associated a function $Y_\lambda(x)$
to a Young diagram $Y_\lambda$ according to \eqref{def-Y}.  

\begin{definition}
We say that  a function $y:\BR_+\to\BR_+$ such that $\int_0^\infty y(x)dx=1$ is the {\em limit shape} of random Young diagrams governed by 
probability measures $M_n$ on $\CP_n$ with the scaling sequence $\beta_n$ if for any $x_0>0$ and $\eps>0$
\begin{equation}\label{ls}
\lim_{n\to\infty}\sup_{x>x_0}M_n\bigl\{\lambda\in\CP_n:\bigl|\tfrac{\beta_n}{n}Y_\lambda(x\beta_n)-y(x)\bigr|>\eps \bigr\}=0.
\end{equation} 
\end{definition}

In words, the existence of the limit shape $y(x)$  means that 
the scaled Young diagrams of partitions of $n$ (by $\beta_n$ along $x$ axis and by $n/\beta_n$ along $y$ axis, so 
that the scaled Young diagram has unit area) become close to a nonrandom figure $\{(x,y):x\ge 0$ and $0\le y\le y(x)\}$
with high $M_n$-probability, as $n$ grows to infinity.  This can be viewed as an analog of the law of large numbers for 
random partitions.

It is well known that random partitions governed by the uniform measures $M_n^{\text{uni}}$ 
have the limit shape with the scaling sequence $\beta_n=\sqrt{n}$, which leads to the symmetric scaling by $\sqrt{n}$ along
both axes.  The limit shape is defined by the equation
\begin{equation}\label{ls-uni}
\ee^{-\pi x/\sqrt{6}}+\ee^{-\pi y(x)/\sqrt{6}}=1.
\end{equation}
The existence of the limit shape has been shown by Temperley \cite{Temperley} using physical argumentation and then
rediscovered by Vershik  \cite{Vershik-Kerov77,V-FAA}.  Now many proofs of this fact are known, including
a purely combinatorial proof by Petrov~\cite{Petrov}.

We do not provide a new rigorous proof here but explain heuristically why the limit shape is formed in this model and
has the form \eqref{ls-uni}.   Note that, in fact, our framework allows to speak about the \textit{strong limit shape},
that is the claim that for $y(x)$ defined by \eqref{ls-uni}
\begin{equation}\label{sls}
\lim_{j\to\infty} \frac{1}{\sqrt{N(\hLa_j)}}Y_{\hLa_j}\Bigl(x\sqrt{N(\hLa_j)} \Bigr)= y(x)\quad\text{a.s.}
\end{equation}
which we believe to be true, maybe even uniformly over $x>0$, however do not prove.
The main observation is that, according to Theorem~\ref{th:hat-Lambda},
the distribution of a new rectangle $ K_j\times R_j$ added to $\hLa_j$ 
is stationary (does not depend on $j$) and
depends on $\hLa_j$ only through
its weight $N(\hLa_j)$.  Moreover, as $n=N(\hLa_j)$ becomes large, the joint distribution 
of $(K_j/\sqrt{n},R_j)$ converges weakly to the distribution described in Proposition~\ref{prop:kr-limits}.

First we check that the jumps of the chain $N(\hLa_j)$ satisfy
$N(\hLa_{j+1})-N(\hLa_j)=O(N(\hLa_j)^{1/2+\eps})$ a.s.\ as $j\to\infty$, for any $\eps>0$.

\begin{lemma}\label{lem:area-jump}
For any $\eps>0$,  there exists an a.s.\ finite random step $J_0=J_0(\eps)$ of the Markov chain
$(N(\hLa_j),K_j,R_j)_{j\ge 0}$, such that 
$K_jR_j\le N(\hLa_j)^{1/2+\eps}$ for all $j\ge J_0$.
\end{lemma}

\begin{proof}
From the second claim of Proposition \ref{prop:kr-limits}(iii), for any $\nu$ and $\eps>0$, applying Markov's inequality 
to $(K_jR_j)^\nu$ we see that for some constant $c=c(\nu,\eps)$
the bound  $\BP[K_jR_j\ge n^{1/2+\eps}|N(\hLa_j)=n]\le c/n^{\nu\eps }$ holds for large $n$.  Since 
$N(\hLa_j)\ge j$ for any $j\ge 0$ by construction, this yields
 $\BP[K_jR_j\ge n^{1/2+\eps}|N(\hLa_j)=n]\le c/j^{\nu\eps }$. Applying this inequality for $\nu>1/\eps$ and using the 
Borel--Cantelli lemma, we see that just finitely many of events $\{K_jR_j\ge N(\hLa_j)^{1/2+\eps}\}$
occur in a realization of the Markov chain $(N(\hLa_j),K_j,R_j)_{j\ge 0}$, with probability~1. So one
can take $J_0-1$ to be the step on which the last such event occurs.
\end{proof}

Denote for short $N_j=N(\hLa_j)$.  
The new diagram $Y_{\hLa_{j+1}}$ is obtained from $Y_{\hLa_j}$ by insertion of a rectangle $K_j\times R_j$, which
can be written as 
\begin{equation}\label{Y-evo}
Y_{\hLa_{j+1}}(x)=Y_{\hLa_j}(x)+R_j\boldone_{\{K_j\ge x\}}, \qquad x\ge 0.
\end{equation}
Different scalings should be applied to $Y_{\hLa_{j}}$ and $Y_{\hLa_{j+1}}$ to obtain diagrams of unit area: 
the former is scaled by $\sqrt{N_{j}}$, while the latter is scaled by $\sqrt{N_{j+1}}$.  However from Lemma~\ref{lem:area-jump}
we know that not only $N_j\sim N_{j+1}$ a.s.\ as $j\to\infty$, but moreover
$N_j\sim N_{j+m_j}$ a.s.\ as long as $m_j=O(N_j^{1/2-\eps})$ as $j\to\infty$, for any $\eps>0$.  On the other hand,
Proposition \ref{prop:kr-limits} allows to conclude that while $N_j\sim N_{j+m_j}$ the distribution of 
$\bigl(K_{j+m_j}/\sqrt{N_{j+m_j}},R_{j+m_j}\bigr)$ is close to the distribution of $\bigl(K_{j}/\sqrt{N_{j}},R_j\bigr)$ for large $j$, and close 
to their joint limit described in Proposition \ref{prop:kr-limits}(i--ii).

Item (iv) of Proposition \ref{prop:kr-limits} ensures that $\BE\bigl[R_i\boldone_{\{K_i\ge x\sqrt{N_i}\}}\bigr]\to h(x)$
as $j\to\infty$, where $h(x)$ is defined in \eqref{h-def}, for all $i=j,j+1,\dots,j+m_j$ such that   $m_j=O(N_j^{1/2-\eps})$, for any $\eps>0$.
Hence from the strong law of large numbers we can conclude that adding several random summands $R_i\boldone_{\{K_i\ge x\sqrt{N_i}\}}$,
for $i=j,j+1,\dots,j+m_j$,
to a scaled Young diagram $\frac{1}{\sqrt{N_j}}Y_{\hLa_j}$ is approximately the same as 
adding a non-random quantity $\frac{m_j}{\sqrt{N_j}}\BE \bigl[R_j\boldone_{\{K_j\ge x\sqrt{N_j}\}}\bigr]$, as long as $j$ is large and 
$m_j$ is also large but smaller than $N_j^{1/2-\eps}$.  

These considerations show that as $j$ becomes large, the random evolution
of scaled Young diagram $y_j(x):=\frac{1}{\sqrt{N_j}}Y_{\hLa_j}(x\sqrt{N_j})$ becomes better and better approximated by the following 
deterministic evolution.  For a real number $r>0$ consider the operator $A_r$ which acts on functions $\BR_+\to\BR_+$
with unit integral:
\begin{equation}\label{A_r}
A_r y(x)=\frac{1}{\sqrt{1+r\eta}}\bigl(y(x\sqrt{1+r\eta})+r h(x\sqrt{1+r\eta})\bigr),\qquad
\eta=\int_0^\infty h(x)dx=\frac{2\sqrt{6}}{\pi},
\end{equation}
that is it adds $rh(x)$ and then rescales both axes by the same factor so that the integral of the resulting function is again 1.
For $r$ close to zero this evolution can be approximated as
\[
A_r y(x) = y(x)+\bigl(h(x)+\tfrac{\eta}{2} x y'(x)-\tfrac{\eta}{2} y(x)\bigr)r+O(r^2), \qquad r\downarrow 0.
\]
The dynamics $y_j(x)\mapsto A_r y_j(x)$ approximates the evolution of $y_j(x)$ for small $r=\smash[b]{\frac{m_j}{\sqrt{N_j}}}$.  
This implies that the
fixed point $y$ of this dynamics should satisfy 
\begin{equation}\label{y-diff-eq}
h(x)+\tfrac{\eta}{2} x y'(x)-\tfrac{\eta}{2} y(x)=0.
\end{equation}
For $h(x)$ given by \eqref{h-def}, it is easy to see that the only solution of \eqref{y-diff-eq} vanishing
at $\infty$ is 
\[
y(x)=-\tfrac{\sqrt{6}}{\pi}\log(1-\ee^{-\pi x/\sqrt{6}}), \qquad x>0.
\]
This is Vershik's limit shape \eqref{ls-uni}.

Similar replacement of a random growth procedure with an approximating it deterministic procedure 
is instrumental in \cite{Kerov99}, but the reasoning there is different from ours.

\section{Adjacent level growth}\label{sec:adding-a-cell}
The major drawback of the growth model defined and investigated in the previous sections is that the growing diagram 
$\hLa_j$ jumps over many levels of the Young graph $\CP$.  A random walk on $\CP$ which visits all levels and
reaches level $n$ on step $n$, that is for which a new partition $\la\vdash n+1$ is obtained from $\mu\vdash n$
by increasing one of parts by 1, or adding a new part 1, would be much easier to analyze. 
Actually, the question of existence of such a random growth procedure is reduced to a study
of adjacent level $\CP_n$ and $\CP_{n+1}$ of the Young graph $\CP$, for each $n\ge 0$.

If $n$ is fixed, we are interested in finding the probabilities $\tp_{\mu,\la}$
of transition from $\mu$ to $\la$, for each $\mu\vdash n$ and $\la\vdash n+1$ satisfying $\mu\nearrow\la$,
such that taking $\mu$ uniformly at random on $\CP_n$ and moving it to $\la\nwarrow\mu$ with probability $\tp_{\mu,\la}$
gives $\la$ with the uniform distribution on $\CP_{n+1}$.   
This is equivalent to finding a non-negative solution to a linear system of equations, which
becomes symmetric in variables $\tq_{\mu,\la}=p(n+1)\tp_{\mu,\la}$:
\begin{alignat}{2}\label{q-system-a}
\sum_{\la:\la\nwarrow\mu}\tq_{\mu,\la}&=p(n+1),\qquad &&\forall \mu\vdash n,\\
\label{q-system-b}
\sum_{\mu:\mu\nearrow\la}\tq_{\mu,\la}&=p(n),\qquad &&\forall \la\vdash n+1.
\end{alignat} 
Equation \eqref{q-system-a} claims that $(\tp_{\mu,\la})_{\la:\la\nwarrow\mu}$ is a probability distribution for each $\mu\vdash n$,
while equation \eqref{q-system-b} ensures that taking $\mu$ uniformly at random from $\CP_n$
and making a step to $\la\nwarrow\mu$ with probability $\tp_{\mu,\la}$ one gets $\la$ 
uniformly distributed on $\CP_{n+1}$.

The necessary and sufficient condition for the system \eqref{q-system-a}--\eqref{q-system-b}, for any bipartite graph 
with parts $\CP_n$ and $\CP_{n+1}$ of cardinalities $p(n)$ and $p(n+1)$, correspondingly, to 
have a non-negative solution is given by a variant of the so-called Supply-Demand Theorem of linear programming,
see, e.g., \cite[Cor.\ 2.1.5]{Lovasz-Plummer}.  This theorem is an easy consequence of the Max-Flow Min-Cut
Theorem \cite[Th.\ 2.1.4]{Lovasz-Plummer} or of Strassen's theorem on stochastic domination \cite{Strassen}.
For a set of partitions $A\subset\CP$ denote its size by $|A|$  and let
\[
A^\uparrow \eqdef \{\la\in\CP:\la\nwarrow\mu\text{ for some }\mu\in A\}.
\]

\begin{lemma}[Uniform version of the Supply-Demand Theorem]
The system \eqref{q-system-a}--\eqref{q-system-b} has a non-negative solution
if and only if for any $A\subset \CP_n$
\begin{equation}\label{density}
\frac{|A^\uparrow|}{|A|}\ge \frac{p(n+1)}{p(n)}.
\end{equation}
\end{lemma}

\begin{remark*}
An equivalent form in terms of uniform measures of the inequality \eqref{density} is 
\begin{equation}\label{Mn-increase}
 M_{n+1}^{\text{\textup{uni}}}(A^\uparrow)\ge M_n^{\text{\textup{uni}}}(A).
\end{equation} 
\end{remark*}

We were unable to check whether the condition \eqref{density} holds for any two neighboring levels of the Young graph.
Note that in the related problem for strict partitions of an integer $n$, which are defined as decreasing sequences
$(\la_1,\dots,\la_\ell)$ of positive integers that sum up to $n$, a one-box growth procedure leading 
to uniform measures on levels can not be constructed for all neighboring levels.  Indeed, if $n=k(k+1)/2$ for some $k\in\BN$, there exists a
``triangular'' strict partition $\Delta=(k,k-1,\dots,1)$ for which just one box can be added increasing the maximal part $k$.
Since the number of strict partitions of $n$ strictly increase in $n\ge 4$, taking $A=\{\Delta\}$ one gets
$|A^\uparrow|/|A|=1$ which violates \eqref{density}.

\section{Equivalence of ensembles}\label{sec:equiv}
Let $X:\CP\to \BR$ be some functional on the set of all partitions.  Taking $\la$ at random with some distribution,
we obtain a random variable $X$.  We shall study functionals which do not increase on trajectories $(\la(t))_{t>0}$.

In the Introduction we have defined the inclusion order $\subseteq$ on partitions (inherited from the corresponding Young diagrams).
The growth procedure defined in Section~\ref{sec:gce} induces a weaker order $\preccurlyeq$:  $\mu\preccurlyeq\la$\label{preccurlyeq} if and only if
$\BP[\la(t_1)=\la,\la(t_2)=\mu]>0$ for some $t_1\le t_2$ (so that $\la(t)$ decreases a.s.\ in this order as $t$ grows). 
In other words, $\mu\prec\la$ iff one can obtain $\la $ from $\mu$ adding further parts to $\mu$ 
(and not changing or removing existing parts). For instance, $(2,1)\subset(2,2)$ but $(2,1)\nprec(2,2)$. Recall that a subset 
$\CB$ of a partially ordered set $(\CP,\preccurlyeq)$ is an \textit{upper order ideal} if from $\mu\in\CB$ and $\mu\preccurlyeq\la$ if follows that 
$\la\in\CB$.  
We observe the following monotonicity property of upper order ideals in $(\CP,\preccurlyeq)$ which will serve as a replacement for \eqref{Mn-increase}.

\begin{lemma}\label{lem:Mn-monotonicity}
Let $\CB\subset \CP$ be an upper order ideal 
of $(\CP,\preccurlyeq)$.   Then, for each $n\ge 1$,
\begin{equation}\label{Mn-monotonicity}
\sum_{m=0}^{n-1} \frac{p(m)\sigma_1(n-m)}{np(n)} M_{m}^{\text{\textup{uni}}}(\CB)
\le M_{n}^{\text{\textup{uni}}}(\CB)\le \sum_{m=n+1}^\infty \frac{g(m)\sigma_1(m-n)}{m g(n)}M_{m}^{\text{\textup{uni}}}(\CB),
\end{equation}
where $p(n)=|\CP_n|$ is the number of partitions of~$n$, $g(n)$ is the visit probability defined in \eqref{gn}, \eqref{g-zero-n}, and  $\sigma_1(h)$ is the sum of divisors function.
\end{lemma}

\begin{remark*}
According to \eqref{npn-recursion} and \eqref{gn-recurstion}, for $\CB=\CP$ both inequalities in \eqref{Mn-monotonicity}  become
equalities. Hence \eqref{Mn-monotonicity} claims that $M_{n}^{\text{\textup{uni}}}(\CB)$ lies between certain convex combination of $M_{m}^{\text{\textup{uni}}}(\CB)$, with $m<n$ and with $m>n$.
\end{remark*}

\begin{proof}
Since for any $n\ge 0$ and $\lambda\in\CP_n$, $\BP[\exists j:\hLa_j=\la]=g(n)$ and $\BP[\exists j:\hLa_j\in\CP_n]=p(n)g(n)$, we can write,
for $n\ge 1$, 
\begin{align}\notag 
M_n^{\text{uni}}(\CB)&=\frac{1}{g(n)p(n)}
\BP[\exists j:\hLa_j\in\CB, N(\hLa_j)=n]\\ \notag 
&=\frac{1}{g(n)p(n)}\sum_{m=0}^{n-1}\BP[\exists j:\hLa_j\in\CB, N(\hLa_j)=n, N(\hLa_{j-1})=m]\\ \notag 
&\ge\frac{1}{g(n)p(n)}\sum_{m=0}^{n-1}\BP[\exists j:\hLa_j\in\CB, N(\hLa_j)=n, N(\hLa_{j-1})=m,\hLa_{j-1}\in\CB]\\ \label{Mn-ineq}
&=\frac{1}{g(n)p(n)}\sum_{m=0}^{n-1}\BP[\exists j:N(\hLa_j)=n, N(\hLa_{j-1})=m,\hLa_{j-1}\in\CB].
\end{align}
The last equality holds because $\CB$ is an upper order ideal, and $\hLa_j\in\CB$ implies $\hLa_{j'}\in\CB$ for all $j'\ge j$.
Using Markov's property and the fact that the distribution of a jump $N(\hLa_j)-N(\hLa_{j-1})$ depends on $\hLa_{j-1}$ only
through its weight $N(\hLa_{j-1})$, we obtain
\[\BP[\exists j:N(\hLa_j)=n, N(\hLa_{j-1})=m,\hLa_{j-1}\in\CB]=\BP[\exists j:N(\hLa_{j-1})=m,\hLa_{j-1}\in\CB]\,\hq_{m,n},
\]
where $\hq_{m,n}=\frac{g(n)\sigma_1(n-m)}{ng(m)}$ is defined in \eqref{hq}.  Plugging this into \eqref{Mn-ineq} and using
the representation in the first line of \eqref{Mn-ineq} (with $m$ instead of $n$) we get the lower bound in~\eqref{Mn-monotonicity}.

To get the upper bound we consider a (non-standard) Markov chain $\Lambda_j=\hLa_{-j}$\label{Lambda-j}, indexed by non-positive integers $j$,
which has transition probabilities $\BP[\Lambda_{j+1}=\la|\Lambda_j=\mu]=p_{\la,\mu}$.  It has no initial distribution (at
``time'' $j=-\infty$), but has terminal distribution $\Lambda_0=\varnothing$ a.s.  An alternative approach could 
be to consider chains $(\Lambda_j^\tau)_{j=0,\dots,J(\tau)}$ and then pass to a limit $\tau\downarrow 0$, but we 
find this approach more involved.

Using the chain $(\Lambda_j)_{j\le 0}$ we can proceed similarly to \eqref{Mn-ineq}.  We have
\begin{align*}
M_n^{\text{uni}}(\CB)&=\frac{1}{g(n)p(n)}\BP[\exists j:\Lambda_j\in\CB, N(\Lambda_j)=n]\\
&=\frac{1}{g(n)p(n)}\sum_{m=n+1}^{\infty}\BP[\exists j:\Lambda_j\in\CB, N(\Lambda_j)=n,\Lambda_{j-1}\in\CB,N(\Lambda_{j-1})=m]\\
&\le \frac{1}{g(n)p(n)}\sum_{m=n+1}^{\infty}\BP[\exists j:N(\Lambda_j)=n,\Lambda_{j-1}\in\CB,N(\Lambda_{j-1})=m]\\
&=\frac{1}{g(n)p(n)}\sum_{m=n+1}^{\infty}g(m)p(m)M_m^{\text{uni}}(\CB)\,q_{m,n}. 
\end{align*} 
Here $q_{m,n}=\frac{p(n)\sigma_1(m-n)}{mp(m)}$, see \eqref{qnm}.  In the second line we have included
the condition $\Lambda_{j-1}\in\CB$ because it follows from $\Lambda_j\in\CB$.  The upper bound in \eqref{Mn-monotonicity} 
is obtained after cancellation.
\end{proof}

We formulate the following theorem, which is one of our main results, in terms of growing 
partition process $\tla(s)$ constructed from the Poisson point process $\tPi$.  Since the
random partition $\tla(s)$ has Fristedt's distribution $\mgc{\exp(-F^{-1}(s))}$, 
the theorem can be easily reformulated in terms of the convergence of random partitions 
governed by $\mgc{q}$, as $q\to1$.  This may seem more traditional reformulation but 
the condition \eqref{a-slow} 
becomes more cumbersome.  Recall that a function $a:\BR_+\to\BR$ is called regularly
varying (at infinity) if $a(cs)/a(s)\to c^r$ as $s\to\infty$, for any $c>0$ and some constant $r\in\BR$.

\begin{thm}\label{th:equivalence}
Let $X(\la)=(X_1(\la),\dots,X_d(\la))$ be a vector such that each $X_i(\la)$ does not decrease 
as a functional on the poset $(\CP,\preccurlyeq)$.  Suppose that 
there exist regularly varying functions $\ta_i(s)$ and $\tb_i(s)>0$ \textup{($i=1,\dots,d$)} such that 
a rescaled random vector 
\begin{equation}\label{Z(s)-def}
Z(s)
\eqdef\bigl(\tb_1(s)(X_1(\tla(s))-\ta_1(s)),\dots,\tb_d(s)(X_d(\tla(s))-\ta_d(s))\bigr)
\end{equation}
converges in distribution, as $s\to\infty$, 
to a random vector $Z\in \BR^d$.  Assume additionally that, as $s\to\infty$, for some function $u(s)=o(s^{1/4})$ but $u(s)\to\infty $,
and for every $i=1,\dots,d$
\begin{gather}\label{a-slow}
\tb_i(s)\bigl(\ta_i(s)-\ta_i(s-u(s)\sqrt{s})\bigr)\to 0.
\end{gather}
Then, if $\Lambda(n)$ is chosen at random from $\CP_n$ with the uniform distribution $M_n^{\text{\textup{uni}}}$
and $s_n$ is defined in \eqref{sn}, the random vector 
\[
Z_n:=\bigl(\tb_1(s_n)(X_1(\Lambda(n))-\ta_1(s_n)),\dots,\tb_d(s_n)(X_d(\Lambda(n))-\ta_d(s_n))\bigr)
\]
converges in distribution, as $n\to\infty$,  to $Z$.
\end{thm}

Before proceeding with the proof let us formulate some results we need.  Lemma~\ref{lem:N-asymp-gc}
was formulated for a sequence Fristedt's grand canonical measures $\mgc{q_n}$, where $n\to\infty$ and $q_n\uparrow 1$ 
is given by \eqref{tn-qn}.  
We reformulate it
for continuously growing parameter $q\uparrow1$; moreover, it is more convenient for us to use a
variable $s\to\infty$ which is connected to $q$ by the one-to-one correspondence $q=\ee^{-t}$; $t=F^{-1}(s)$, see Section~\ref{sec:Ppp}.

\begin{lemma}[Local limit theorem]\label{lem:llt} As $s\to\infty$, the following asymptotics holds:
\begin{align}
\label{N-tla-mean}
\BE [N(\tla(s))]&=\frac{6s^2}{\pi^2}+O(s\log s),\\ \label{N-tla-var}
\operatorname{Var}[N(\tla(s))]&=\frac{72 s^3}{\pi^4}+O(s^2\log s).
\end{align}
Let $n(s)$ be an integer-valued function of $s\in[0,\infty)$ such that 
\[
\lim_{s\to\infty} \frac{n(s)-\BE[N(\tla(s))]}{\sqrt{\operatorname{Var}[N(\tla(s))]}}=w.
\]
Then, uniformly in $|w|<W$ for any fixed constant $W$,
\begin{equation}\label{llt}
\lim_{s\to\infty}\sqrt{\operatorname{Var}[N(\tla(s))]}\,\BP\bigl[N(\tla(s))=n(s)\bigr]=\frac{1}{\sqrt{2\pi}}\ee^{-\frac{w^2}{2}}.
\end{equation}
\end{lemma}
\noindent 
Lemma~\ref{lem:llt} is a specialization of \cite[Lemma 10]{Yakubovich12}.  See also \cite{Baez-Duarte}.
Uniformity is not claimed in \cite{Yakubovich12} but can be easily checked by tracing the proof. Then the asymptotics 
for $t\downarrow 0$ should be rewritten in terms of $s\to\infty$ using \eqref{f-F-s-asymptotics}.  The details
are omitted.\qed

\begin{remark*}
The asymptotic relation \eqref{mgc-N=n} of Lemma~\ref{lem:N-asymp-gc} is a partial case of \eqref{llt} for $w=0$,
calculated for $s=s_n$ as defined in~\eqref{sn}.
\end{remark*}

\begin{corollary}\label{cor:Chebyshev}
If $u(s)\to\infty $ as $s\to\infty$, then
\begin{equation}\label{Chebyshev}
\BP\bigl[\bigl|N(\tla(s))-6s^2/\pi^2\bigr|>u(s)s^{3/2}\bigr]=O(u(s)^{-2})\to 0.
\end{equation}
\end{corollary}
\noindent
This is just an application of the Chebyshev inequality to $N(\tla(s))$.\qed

\smallskip

The coefficients of the right convex combinations in \eqref{Mn-monotonicity} are 
jump probabilities from any partition of $n$ to some partition of $m>n$ for the backward (increasing)
Markov chain $(\hLa_j)_{j\ge 0}$, see Proposition~\ref{prop:hN-increasing}.  The left convex combination
coefficients can be interpreted as jump probabilities from a uniform random partition of $n$ to
some partition of $m<n$ for the forward (decreasing) Markov chain $(\Lambda_j)_{j\le 0}$. In both cases,
for large $n$ 
the major contribution to these convex combinations give $m$ such that $|n-m|$ has order $\sqrt{n}$. 
This allows to bound tail sums for these coefficients.

\begin{lemma}\label{lem:convex-combinations}
For any $\alpha\in(1/2,1)$ and any $v>0$, as $n\to\infty$
\[
\sum_{h>n^\alpha}\frac{p(n-h)\sigma_1(h)}{np(n)}=O(n^{-v}), \qquad
\sum_{h>n^\alpha}\frac{g(n+h)\sigma_1(h)}{(n+h) g(n)}=O(n^{-v}).
\]
\end{lemma}

\begin{proof}
The second sum is exactly the probability $\BP[K_jR_j>n^\alpha|N(\hLa_j)=n]$ 
that the backward chain jumps from any partitions $\hLa_j$ of $n$ to a level bigger
than $n+n^\alpha$.  As in the proof of Lemma~\ref{lem:area-jump}, we can apply
Markov's inequality to $(K_jR_j)^m$, where $m>\delta/(\alpha-1/2)$, and use \eqref{EKRm}.
To obtain the first bound without calculation, we can argue that when the backward
chain $\hLa$ jumps from $n$ to $m$, the forward chain $\Lambda$ jumps from $m$ to $n$,
and as long as $m\sim n$ the distributions of these jumps have the same tail behavior.  
\end{proof}

We also need the following claim about level visit probabilities.  Recall that
the chain $\hLa_j$ visits level $\CP_n$ with probability $\gamma(n)\sim\pi/(2\sqrt{6n})$ as 
$n\to\infty$, see Corollary~\ref{cor:level-visit-prob}.

\begin{lemma}\label{lem:visit-two-levels}
If $n,m\to\infty$ and $m-n\ge m^{3/4}r(m)$, where $r(m)\to\infty$ and do not decrease, then the events
$\{N(\hLa_j)=n$ for some $j\}$ and $\{N(\hLa_{j'})=m$ for some $j'\}$ are asymptotically independent:
\begin{equation}\label{gammanm}
\gamma(n,m)\eqdef\BP [\exists j,j': N(\hLa_j)=n,N(\hLa_{j'})=m]\sim \gamma(n)\gamma(m)\sim\frac{\pi^2}{24\sqrt{nm}}.
\end{equation}
\end{lemma}

\begin{proof}
Let us consider the Poisson point process description $\tla(s)$ of the growing partition process.
Then $\gamma(n,m)=\BP[\exists s',s'': N(\tla(s'))=n, N(\tla(s''))=m]$.  
Let $n^+=n^+(m)\eqdef m-m^{3/4}r(m)$; then $n^+\ge n$ by assumption. 
Consider $s_{n^+}$ and  $s_m$ defined by \eqref{sn} and let
\[
\sbar = \frac{s_{n^+}+s_m}{2}.
\]
Then it is easy to check using \eqref{sn}, \eqref{N-tla-mean} and \eqref{N-tla-var} that for all sufficiently large~$n,m$
such that $n\le m-m^{3/4}r(m)$  
inequalities
\[
n+\tfrac13m^{3/4}r(m) \le \BE[\tla(\sbar)]\pm \tfrac13\sqrt{\operatorname{Var}[N(\tla(\sbar))]}r(m)\le m-\tfrac13m^{3/4}r(m)
\]
hold.   Let $\ell$ be such that $|\ell-\BE[\tla(\sbar)]|<\tfrac13\sqrt{\operatorname{Var}[N(\tla(\sbar))]}r(m)$.
Then, given $N(\tla(\sbar))=\ell$, jumps of $(N(\tla(s)))_{s\ge \sbar}$ and $(N(\tla(s)))_{s\le\sbar}$ are conditionally
independent.  Hence the events $\{\exists s'<\sbar: N(\tla(s'))=n\}$ and $\{\exists s''>\sbar: N(\tla(s''))=m\}$
are conditionally independent given $N(\tla(\sbar))=\ell$ for all such $\ell$.

We claim that both conditional probabilities  $\BP[\exists s'<\sbar: N(\tla(s'))=n|N(\tla(\sbar))=\ell]\allowbreak\sim\allowbreak g(n)$ 
and $\BP[\exists s''>\sbar: N(\tla(s''))=m|N(\tla(\sbar))=\ell]\sim g(m)$ as $n\to\infty$ and $m,\ell$ are as above. We shall prove
just the latter claim which is enough for our purposes.  
To show it, we can, following \cite{Gnedin-Yak}, construct two renewal processes 
$(\hN_j^+)_{j\ge j_0}$ and $(\hN_j^-)_{j\ge j_0}$ which start from $\ell=\hN_{j_0}^+=\hN_{j_0}^-$ and move faster and slower than $(N(\hLa_j))_{j\ge j_0}$ 
given $N(\hLa_{j_0})=\ell$, respectively, but still 
have asymptotically equivalent visit probabilities of the state $m$. 

The construction of these stochastic bounds is easy due to the fact that $\hq_{n,n+h}$ eventually decreases in $n$ for any fixed $h$. 
Indeed, after a little algebra we obtain
\begin{equation}\label{ddn-hq}
\frac{d}{dn}\hq_{n,n+h}=
\frac{\ee^{A_n(0)-A_n(h)} \sigma_1(h) \left(u_1(n,h)+u_2(n,h)\right)}{n^2 \sqrt{24 n-1} (24 h+24 n-1)^{3/2}}
\end{equation}
where
\[
u_1(n,h)=24 h \left(12 n \left(A_n(0)-1\right)+1\right),\quad u_2(n,h)=(24 n-1) \left(12 n \left(A_n(0)-A_n(h)-2\right)+1\right).
\]
It is easy to see that $u_1(n,h)>0$, $u_2(n,h)<0$ and $|u_2(n,h)|>u_1(n,h)$ at least for $n>h$, hence the derivative \eqref{ddn-hq} is negative.
Let $(X_j^+)$ and $(X_j^-)$ be two sequences of i.i.d.\ random variables such that $\BP[X_j^+=h]=\hq_{m,m+h}$ and $\BP[X_j^-=h]=\hq_{\ell,\ell+h}$.
Then the following inequalities hold:
\begin{equation}\label{tail-inequalities}
\BP[X_j^+\ge h]\ge \BP[N(\hLa_{j+1})-N(\hLa_j)|N(\hLa_j)=k]\ge \BP[X_j^-\ge h],\qquad h<\ell\le k\le m.
\end{equation}
Hence we can couple triples $(X_j^+,X_j^-,N(\hLa_{j+1})-N(\hLa_j))$ in such a way that $X_j^+\ge N(\hLa_{j+1})-N(\hLa_j)\ge X_j^-$
as long as $\ell\le N(\hLa_j)\le m$. (Here we use Lemma~\ref{lem:convex-combinations} to 
overcome a situation when some of these jumps is greater than $\ell$ for which we have not checked inequalities \eqref{tail-inequalities}:
the probability that it happens for at least one of the considered random variables is smaller than any negative power of $\ell$.)
Hence, given $N(\hLa_{j_0})=\ell$ for some $j_0$, the renewal process $\hN_j^+\eqdef\ell+X_{j_0}^+\dots+X_{j-1}^+$ 
jumps faster than the Markov chain $N(\hLa_j)$ until the latter overcomes the level $m$, and this Markov chain jumps faster than the renewal process $\hN_j^-\eqdef\ell+X_{j_0}^-\dots+X_{j-1}^-$.

According to Proposition~\ref{prop:kr-limits}(iii), the expectations of jumps $X_j^\pm$ for both renewal processes
have the same asymptotics $2\sqrt{6m}/\pi$ as $m\to\infty$, because $\ell\sim m$.  By Lemma~\ref{lem:convex-combinations} and 
assumptions that $m-\ell\ge \tfrac13m^{3/4}r(m)$, each chain with probability close to 1 needs at least $m^{1/4-\varepsilon}$ steps to 
reach level $m$ starting from level $\ell$, for any $\eps>0$ as long as $m$ is large enough.  Since 
$\BP[X_j^\pm=h]>0$ for all $h\in\BN$, by a standard result of
renewal theory the probability that any of $(\hN_j^\pm)$ visits $m$ is asymptotically $\frac{\pi}{2\sqrt{6m}}$.
Since $(N(\hLa_j))$ visits at least as many levels as $(\hN_j^+)$ but no more levels than $(\hN_j^-)$, it follows that 
$\BP[\exists j:N(\hLa_j)=m|N(\hLa_{j_0})=\ell]\sim \frac{\pi}{2\sqrt{6m}}$ as $m\to\infty$.

To finish the proof it remains to apply the total probability formula:
\begin{align*}
\BP[&\exists j',j'':N(\hLa_{j'})=n,N(\hLa_{j''})=m]\\
&=\sum_{\ell= \ell^-}^{\ell^+}\BP[\exists j',j'':N(\hLa_{j'})=n,N(\hLa_{j''})=m|N(\tla(\sbar))=\ell]\,\BP[N(\tla(\sbar))=\ell]\\
&\qquad+\BP[\exists j',j'':N(\hLa_{j'})=n,N(\hLa_{j''})=m|N(\tla(\sbar))<\ell_-]\,\BP[N(\tla(\sbar))<\ell^-]\\
&\qquad+\BP[\exists j',j'':N(\hLa_{j'})=n,N(\hLa_{j''})=m|N(\tla(\sbar))>\ell^+]\,\BP[N(\tla(\sbar))>\ell^+]
\end{align*}
where $\ell^\pm=\BE[N(\tla(\sbar))]\pm \frac13 \operatorname{Var}[N(\tla(\sbar))]r(m)$.
According to Corollary~\ref{cor:Chebyshev} we may expect that the sum gives the main contribution, 
and using the conditional independence of the past and the future
and asymptotic independence of events $\{N(\tla(\sbar))=\ell\}$ and  $\{\exists j:N(\hLa_j)=m\}$ established
above we conclude that the sum is asymptotic to $\gamma(n)\gamma(m)$. Two additional summands
include factors $\BP[N(\tla(\sbar))<\ell^-]$ and $\BP[N(\tla(\sbar))>\ell^+]$ which vanish as $m\to\infty$, but
possibly slower than $\gamma(n)\gamma(m)$.  So we still have to show that they give a contribution which can be 
neglected, which can be accomplished by repeating the same argument with different $\sbar$.  The details are omitted.
\end{proof}

Now we have all results needed to prove Theorem~\ref{th:equivalence}, but let us first give an outline 
of the proof. First we show that the condition \eqref{a-slow} implies that we can find a family of
upper order ideals $\CB(s,z)$ of $(\CP,\preccurlyeq)$, parameterized by
$s\ge0$ and $z\in\BR^d$, such that (i) events $X(\tla(s))\in\CB(s,z)$ determine the
distribution of $X(\tla(s))$ and of the limit, and (ii) $\lim_{s\to\infty}\BP[X(\tla(s'))\in\CB(s,z)]$
depends on $z$ but not on $s'$ as long as $|s'-s|\le u(s)\sqrt{s}$. This is guaranteed by the condition \eqref{a-slow}.
Then we take a large $n$ and consider $\CB(s_n,z)$ with $s_n$ given by  \eqref{sn}.
We iterate bounds of Lemma~\ref{lem:Mn-monotonicity} to provide bounds for $M_n^{\text{uni}}(\CB(s_n,z))$ 
as convex combinations of $M_m^{\text{uni}}(\CB(s_n,z))$ with $m$ lying far from $n$. Let us concentrate on 
the upper bound, as the lower bound is considered in details in the proof below. For each $\ell>n$, we
find an upper bound for $M_n^{\text{uni}}(\CB(s_n,z))$ as a convex combination of $M_m^{\text{uni}}(\CB(s_n,z))$ with
$m\ge \ell$.  Then we multiply this inequality by the probability that the weight of accelerated partition growth process
$N(\tla(s_n+u(s_n)\sqrt{s_n}))$ is $\ell$ and sum over $\ell>n$. It turns out that $\BP[N(\tla(s_n+u(s_n)\sqrt{s_n}))=m]$ (i)
behaves regularly due to the local limit theorem (Lemma~\ref{lem:llt}), and (ii) is negligibly small if $m-n$ is not big
enough for application of Lemma~\ref{lem:visit-two-levels}.  This lemma allows us then to change the order of summation and
show that the upper bound is asymptotic to $\BP[X(\tla(s_n+u(s_n)\sqrt{s_n}))\in\CB(s_n,z)]$. But this probability has the 
same $n\to\infty$ limit as $\BP[X(\tla(s_n))\in\CB(s_n,z)]$, as it has been established in the beginning of the proof. The lower bound is obtained
similarly using the delayed growth process $\tla(s_n-u(s_n)\sqrt{s_n})$, as detailed in the proof below.

\begin{proof}[Proof of Theorem~\ref{th:equivalence}]
Let us start with clarifying the meaning of the condition \eqref{a-slow}. Denote for short $s^\pm\eqdef s\pm u(s)\sqrt{s}$.
The condition \eqref{a-slow} together with regularity assumptions on $\ta_i$, $\tb_i$ guarantee that taking
the delayed (that is $\tla(s^-)$) or accelerated ($\tla(s^+)$) random
partition  instead of $\tla(s)$ in the definition \eqref{Z(s)-def} of $Z(s)$ 
implies the same limit $Z$ after rescaling suitable for $\tla(s)$. To check this 
let us introduce some notation.
For a point $z=(z_1,\dots,z_d)\in\BR^d$ and $s\ge0$ denote 
\[
\CB(s,z)\eqdef\{\la\in\CP:(\tb_1(s)(X_1(\la))-\ta_1(s))>z_1,\dots,\tb_d(s)(X_d(\la)-\ta_d(s))>z_d\}\,.
\]
Note that from the monotonicity assumptions on $X_j$ it follows that $\CB(s,z)$ is an upper order ideal
in $(\CP,\preccurlyeq)$. 
Let $G(z)$ be the ``tail'' distribution function of the limiting distribution of~$Z$:
\[
G(z)\eqdef\BP[Z_1>z_1,\dots,Z_d>z_d].
\]
The values $G(z)$ at continuity points $z$ determine the distribution of $Z$.

Let $z$ be a continuity point of $G$; then by assumption $\BP[\tla(s)\in\CB(s,z)]\to G(z)$ as $s\to\infty$.  
Since $u(s)\sqrt{s}=o(s)$, it follows that $s^\pm\to\infty$ and 
we also have $\BP[\tla(s^\pm)\in \CB(s^\pm,z)]\to G(z)$ as $s\to\infty$, but we need a different claim.
Consider 
\[
Z^\pm(s)\eqdef\bigl(\tb_1(s)\bigl(X_1(\tla(s^\pm))-\ta_1(s)\bigr),\dots,\tb_d(s)\bigl(X_d(\tla(s^\pm))-\ta_d(s)\bigr)\bigr),
\]
that is, use the scaling suitable for $\la(s)$ for its accelerated or delayed version $\la(s^\pm)$.  For each $i$, the event
\begin{multline*}
\bigl\{\tb_i(s)\bigl(X_i(\tla(s^\pm)-\ta_i(s)\bigr)>z_i\bigr\}\\
=\bigl\{\tb_i(s^\pm)\bigl(X_i(\tla(s^\pm))-\ta_i(s^\pm)\bigr)
  >\tfrac{\tb_i(s^\pm)}{\tb_i(s)}\bigl(z_i+\tb_i(s)(\ta_i(s)-\ta_i(s^\pm)\bigr)\bigr\}.
\end{multline*} 
Due to the assumed regular variation of $\ta_i$ and $\tb_i$ combined with the condition \eqref{a-slow}, the lower bound in 
the second line converges to $z_i$ as $s\to\infty$, so 
\begin{equation}\label{pm-same-limits}
\BP\bigl[\tla(s^\pm)\in \CB(s,z)\bigr]\to G(z), 
\end{equation}
because $G$ is continuous at~$z$.  

\smallskip 

Our aim is to show that, as $n\to\infty$, the sequence $M_n^{\text{uni}}[\CB(s_n,z)]$ (where $s_n$ is 
defined by \eqref{sn}) has a limit which is equal to $G(z)$.  
Let $n$ be large, and denote for short $\CB=\CB(s_n,z)$.  For each
$\ell<n$ we can iterate the left inequality in \eqref{Mn-monotonicity}, applying it first to $M_{n-1}^{\text{uni}}(\CB)$,
then to $M_{n-2}^{\text{uni}}(\CB)$ and so on up to $M_{\ell+1}^{\text{uni}}(\CB)$, and collecting
terms at $M_{m}^{\text{uni}}(\CB)$ for $m\le \ell$ to obtain
\begin{equation}\label{Mn-lower-itarated}
M_n^{\text{uni}}(\CB)\ge \sum_{m=0}^\ell \gamma_{n,\ell,m}M_m^{\text{uni}}(\CB),\qquad \ell<n.
\end{equation}
More formally, we define recursively
\[
\gamma_{n,n-1,m}\eqdef q_{n,m}\;\, (m\le n-1),\quad\; \gamma_{n,\ell-1,m}\eqdef \gamma_{n,\ell,m}+\gamma_{n,\ell,\ell}\,q_{\ell,m}\;\, (0\le m\le\ell-1<n-1).
\]
Recall the definition of the decreasing chain $(\Lambda_j)_{j\le0}$ made in the proof of Lemma~\ref{lem:Mn-monotonicity}:
$\Lambda_j=\hLa_{-j}$, $j\le 0$.  The value $\gamma_{n,\ell,m}$ can be interpreted in terms of this chain as the conditional probability that
given $(\Lambda_j)_{j\le0}$ visits $\CP_n$, the first level $\CP_k$ with $k\le \ell$ it visits is exactly $\CP_m$.  
Hence, for any $\ell<n$, 
\begin{equation}\label{gammanellm-sum}
\sum_{m=0}^\ell\gamma_{n,\ell,m}=1.
\end{equation}

We can also write
\begin{equation}\label{gammaellnm}
\gamma_{n,\ell,m}=\sum_{k=\ell+1}^{n}\BP[\exists j':N(\Lambda_{j'})=k|N(\Lambda_j)=n]\,q_{k,m}
=\sum_{k=\ell+1}^n\frac{\gamma(k,n)\,q_{k,m}}{\gamma(n)}
\end{equation}
by considering the last level $\CP_k$ (with $k>\ell$) the chain  $(\Lambda_j)_{j\le0}$ visits before it reaches the level $\CP_m$ (with $m\le\ell$),
given it visits $\CP_n$.
Here $\gamma(k,n)$ and $\gamma(n)$ were defined in \eqref{gammanm} and \eqref{visit-prob}, correspondingly.
Suppose now that $n\ge k>\ell>n/2$ and $\ell-m\ge n^{5/8}$. Then, according to Lemma~\ref{lem:convex-combinations} applied with $\alpha=5/8$,
for any $v>0$ we have $\sum_{m<\ell-n^{5/8}}q_{k,m}=O(n^{-v})$, so taking just the last $\lfloor{n^{5/8}}\rfloor+1$ ($\lfloor\cdot\rfloor$ is the floor function) 
summands in the sum \eqref{gammanellm-sum}
we gain an error smaller than any negative power of $n$: 
\begin{equation}\label{gammanellm-sum-restricted}
\sum_{m=\ell-\lfloor n^{5/8}\rfloor}^{\ell}\gamma_{n,\ell,m}=1-O(n^{-v}),\qquad n\to\infty.
\end{equation}
We shall need a similar sum over $\ell$, which is possible to estimate asymptotically when $n,m\to\infty$ but
$n-m$ is sufficiently large.  Using \eqref{gammaellnm} and \eqref{q-hq} we obtain
\[
\sum_{\ell=m}^{m+\lfloor n^{5/8}\rfloor}\!\!\gamma_{n,\ell,m}
=\sum_{\ell=m}^{m+\lfloor n^{5/8}\rfloor}\sum_{k=\ell+1}^n\frac{\gamma(k,n)\,q_{k,m}}{\gamma(n)}
=\sum_{\ell=m}^{m+\lfloor n^{5/8}\rfloor}\sum_{k=\ell+1}^n\frac{\gamma(k,n)\gamma(m)\,\hq_{m,k}}{\gamma(k)\gamma(n)}.
\]
According to Lemma~\ref{lem:convex-combinations}, if $\ell\sim n\to\infty$ then changing the inner summation
to $k$ from $\ell+1$ to $\ell+\lfloor n^{5/8}\rfloor$ will introduce an error not exceeding any negative power on $n$.
If $(n-m)n^{-3/4}\to \infty$, $\ell\le m+n^{5/8}$ and $k\le\ell+n^{5/8}$, then $(n-k)n^{-3/4}\to\infty$ 
and we can apply Lemma~\ref{lem:visit-two-levels} to assert that $\gamma(k,n)\sim \gamma(k)\gamma(n)$.  
Using again Lemma~\ref{lem:convex-combinations} and Proposition~\ref{prop:kr-limits}(iii), we see that as $n> m\to\infty$ 
\[
\sum_{\ell=m}^{m+\lfloor n^{5/8}\rfloor}\sum_{k=\ell+1}^{\ell+\lfloor n^{5/8}\rfloor}\hq_{m,k}
\sim \sum_{\ell=m}^{\infty}\sum_{k=\ell+1}^{\infty}\hq_{m,k} = \BE[N(\hLa_{j+1})-N(\hLa_j)|N(\hLa_j)=m]\sim \frac{1}{\gamma(m)}.
\]
Hence we obtain
\begin{equation}\label{sum-gammanlm}
\lim_{\substack{\vphantom{\lfloor}n,m\to\infty\\ (n-m)n^{-3/4}\to \infty}}\sum_{\ell=m}^{m+\lfloor n^{5/8}\rfloor}\!\!\gamma_{n,\ell,m}=1 .
\end{equation} 

\smallskip 

Let us multiply \eqref{Mn-lower-itarated} by $\BP[N(\tla(s_n^-))=\ell]$ and sum over all $\ell<n$ to obtain
\begin{equation}\label{Mn-lower-bound-weighted}
M_n^{\text{uni}}(\CB)\BP[N(\tla(s_n^-))<n]\ge \sum_{\ell=0}^{n-1}  \sum_{m=0}^\ell \gamma_{n,\ell,m}M_m^{\text{uni}}(\CB)\,\BP[N(\tla(s_n^-))=\ell].
\end{equation}
Note that by Lemma~\ref{lem:llt}, for $s_n$ is defined by \eqref{sn} we obtain, as $n\to\infty$,
\begin{align}\notag
m_n^-\eqdef\BE[N(\tla(s_n^-))]&=\frac{6(s_n^-)^2}{\pi^2}+O(s_n^-\log s_n^-)\\
\label{mn-}
&=n-\frac{2\cdot 6^{1/4}}{\sqrt{\pi}}u(s_n)n^{3/4}
+O\bigl(\sqrt{n}\max\{\log n,u(s_n)^2\}\bigr),\\
\label{sigman-}
(\sigma_n^-)^2\eqdef\operatorname{Var}[N(\tla(s_n^-))]&\sim \frac{2\sqrt{6}}{\pi}n^{3/2}.
\end{align}
Hence, due to the assumption $u(s)=o(s^{1/4})$, the $O(\cdot)$ term in \eqref{mn-} is asymptotically smaller
than two first terms, thus for any fixed $W$ and all large enough $n$ it holds
\begin{equation}\label{lasn-<n}
m_n^-+W\sigma_n^-<n
\end{equation} 
and, moreover, 
\begin{equation}\label{n-lasn->}
n-(m_n^-+W\sigma_n^-)\ge n^{3/4}\sqrt{u(s_n)},
\end{equation} 
so by Corollary~\ref{cor:Chebyshev}
\[
\BP[N(\tla(s_n^-))<n]\to 1,\qquad n\to\infty.
\]
Let $\delta>0$ be fixed and take $W\eqdef \delta^{-1/2}$.  Then, 
restricting the outer sum in \eqref{Mn-lower-bound-weighted} to $\ell\in[m_n^--W\sigma_n^-,m_n^-+W\sigma_n^-]$ we
gain an error not more than $\delta$ if \eqref{lasn-<n} holds for $n$.  Moreover, according to \eqref{gammanellm-sum-restricted} we can take the inner sum 
in \eqref{Mn-lower-bound-weighted} just over $m$ from $\ell-[n^{5/8}]$ up to 
$\ell$ to gain another error smaller than $\delta$ for any sufficiently large $n$.  Hence we obtain 
\begin{equation}\label{Mn-lower-bound-restricted}
M_n^{\text{uni}}(\CB)\BP[N(\tla(s_n^-))<n]\ge \sum_{\ell=m_n^--W\sigma_n^-}^{m_n^-+W\sigma_n^-}  \sum_{\;m=\ell-\lfloor n^{5/8}\rfloor}^\ell \gamma_{n,\ell,m}M_m^{\text{uni}}(\CB)\,\BP[N(\tla(s_n^-))=\ell]-2\delta.
\end{equation}
Now, for any $\ell$ and $m$ over which the above sum is taken, we can apply Lemma~\ref{lem:llt} (taking $s=s_n^-\to\infty$) 
both with $n(s_n^-)=\ell$ and $n(s_n^-)=m$, and since $\ell-m\le n^{5/8}=o(\operatorname{Var}[N(\tla(s_n^-))])$ the
limit \eqref{llt} is the same in both these cases. So, supposing that $n$ is large enough so that in both 
cases  
\[
\Bigl|\sqrt{\operatorname{Var}[N(\tla(s_n^-))]}\,\BP\bigl[N(\tla(s_n^-))=n(s_n^-)\bigr]-\frac{1}{\sqrt{2\pi}}\ee^{-\frac{w^2}{2}}\Bigr| <\delta/(4W)
\]
we see that replacing $\BP[N(\tla(s_n^-))=\ell]$ with $\BP[N(\tla(s_n^-))=m]$ in the sum \eqref{Mn-lower-bound-restricted} 
produces another error not bigger than $\delta$.  Now we want to change the order of summation in~\eqref{Mn-lower-bound-restricted}, and 
replacing the sum over $\ell\in[{m_n^-}-W\sigma_n^-,{m_n^-}-W\sigma_n^-+n^{5/8}]$ and $m\in[\ell,{m_n^-}-W\sigma_n^-]$
with the sum over $\ell\in[{m_n^-}+W\sigma_n^-,{m_n^-}+W\sigma_n^-+n^{5/8}]$ and $m\in[\ell-n^{5/8},{m_n^-}+W\sigma_n^-]$
again we change the value by no more than $\delta$, when $n$ is large.  This way we obtain
\begin{equation}\label{Mn-lower-bound-order-changed}
M_n^{\text{uni}}(\CB)\BP[N(\tla(s_n^-))<n]\ge \sum_{m=m_n^--W\sigma_n^-}^{m_n^-+W\sigma_n^-}  \sum_{\ell=m\vphantom{m_n^-}}^{m+\lfloor n^{5/8}\rfloor} \gamma_{n,\ell,m}M_m^{\text{uni}}(\CB)\,\BP[N(\tla(s_n^-))=m]-4\delta.
\end{equation}
Finally, we use \eqref{sum-gammanlm} to replace the inner sum over $\ell$ with 1, obtaining another error not bigger than $\delta$ when $n$ is large,
which allows us to represent the sum over $m$ as $\BP[\tla(s_n^-)\in\CB]$, up to $\BP\bigl[\bigl|N(\tla(s_n^-))-m_n^-\bigr|\ge W\sigma_n^-\bigr]\le\delta$.   

Repeating with obvious modifications the above argument using the upper bound given by Lemma~\ref{lem:Mn-monotonicity}
and decomposing over values of $N(\tla(s_n^+))$ instead of $N(\tla(s_n^-))$ we obtain the upper bound for $M_n^{\text{uni}}(\CB)\BP[N(\tla(s_n^+))>n]$.
Combining these two bounds we obtain
\[
\frac{\BP[\tla(s_n^-)\in\CB]-6\delta}{1-\delta}\le M_n^{\text{uni}}(\CB)\le \BP[\tla(s_n^+)\in\CB]+6\delta.
\]
According to \eqref{pm-same-limits} both probabilities above converge to $G(z)$ as $n\to\infty$, 
so $G(z)$ should also be the limit of $M_n^{\text{uni}}(\CB)$ because $\delta>0$ is arbitrary.
\end{proof}

\begin{remark*}
The proof of Theorem~\ref{th:equivalence} would be much simpler if \eqref{Mn-increase} could be established.
\end{remark*}

An application of Theorem~\eqref{th:equivalence} to the study of
moments of a uniform random partition $\la$ of $n$, that is to the sum
$\sum_{i=1}^{\ell(\la)}\la_i^p$, is presented in \cite{Yakubovich23}. 
Another application is given in the next Section.

\section{Odd and even parts in a big random partition}\label{sec:oddeven}

We illustrate the strength of Theorem~\ref{th:equivalence} by the following not very unexpected but new, to the
best of our knowledge, result.  Let 
\begin{equation}\label{odd-even}
\ell_o(\la)=\sum_{k\text{ odd}}C_k(\la),\qquad 
\ell_e(\la)=\sum_{k\text{ even}}C_k(\la)
\end{equation}
be the number of odd and even parts in a partition $\la$. So $\ell_o(\la)+\ell_e(\la)=\ell(\la)$, the length of $\la$.  
We use Theorem~\ref{th:equivalence} to 
show that if $\la$ is a random partition of $n$, then $\ell_o(\la)$ and $\ell_e(\lambda)$ are 
asymptotically independent as $n\to\infty$.  As a byproduct we find their limiting distributions which are, 
somewhat surprisingly, different.  This complements the well known result by Erd\H{o}s and Lehner about the 
limit distribution of the length $\ell(\la)$ of the uniform random partition $\la$ of integer $n$:
\begin{equation}\label{EL-Gumbel}
\lim_{n\to\infty}M_n^{\text{uni}}\{\la\in\CP_n:\tfrac{\pi}{\sqrt{6n}}\ell(\la)-\tfrac12\log n+\log\tfrac{\pi}{\sqrt{6}}\le x\}
=\ee^{-\ee^{-x}},\qquad x\in\BR.
\end{equation}
This result was proved in \cite{EL}, perhaps the first paper on random integer partitions, by purely
combinatorial methods.

\begin{thm}\label{th:odd-even} Let $\ell_o(\la)$ and $\ell_e(\la)$ be defined by \eqref{odd-even}. Then
for any $x,y\in\BR$
\begin{equation}\label{elloe-joint}
\lim_{n\to\infty}M_n^{\text{uni}}
\{\la\in\CP_n:\tfrac{\pi}{\sqrt{6n}}\ell_o(\la)-\alpha(n)\le x,
\tfrac{\pi}{\sqrt{6n}}\ell_e(\la)-\alpha(n)\le y\}=\operatorname{Erfc}(\ee^{-x})\ee^{-\ee^{-2y}}
\end{equation} 
where
\[
\operatorname{Erfc}(u)=\frac{2}{\sqrt{\pi}}\int_u^\infty \ee^{-z^2}dz
\]
and 
\begin{equation}\label{alpha}
\alpha(n)=\tfrac14\log n -\tfrac12\log\tfrac{2\pi}{\sqrt{6}}.
\end{equation}

\end{thm}

\begin{proof}
We first show the corresponding limit theorem for $\lambda(t)$ as $t\downarrow0$ and then reformulate it for $\tla(s)$ 
to check the conditions \eqref{a-slow}.  Obviously, $\ell_o(\la(t))$ and $\ell_e(\la(t))$ are independent for any $t>0$, as they 
are defined in terms of different independent processes $\rho_k(t)$. So we can find their limit distributions separately. 

Throughout this proof we assume that $u$ is non-negative real number,
although all equalities involving $u$ below hold for any complex $u$ with $\Re u\ge 0$,
and some of them hold for $u$ from a wider area in the complex plane.
We can write the Laplace transform of $\rho_k(t)$ explicitly as
\[
\BE \,\ee^{-u\rho_k(t)}= \frac{1-\ee^{-kt}}{1-\ee^{-kt-u}}=\frac{1}{1+\frac{\ee^{-kt}}{1-\ee^{-kt}}(1-\ee^{-u})},
\] 
so the cumulant generating function of the centered random variable is
\begin{equation*}
\log \BE\,\ee^{-u(\rho_k(t)-\BE\rho_k(t))} = -\log\Bigl(1+\frac{\ee^{-kt}}{1-\ee^{-kt}}(1-\ee^{-u})\Bigr)-\frac{u\,\ee^{-kt}}{1-\ee^{-kt}}.
\end{equation*} 
This can be rewritten in a L\'evy--Khintchine type form 
\begin{equation}\label{rho-Levy-Khintchine}
\log \BE\,\ee^{-u(\rho_k(t)-\BE\rho_k(t))} = -\sum_{j=1}^{\infty} (1-\ee^{-u j}+u j)\frac{\ee^{-jkt}}{j}
\end{equation}
which should be treated as an integral with respect to the discrete L\'evy measure $\sum_{j=1}^\infty \frac{\ee^{-jkt}}{j}\delta_j$,
with $\delta_j$ meaning the unit point mass at $j$. Plugging $ut$ instead of $u$ in \eqref{rho-Levy-Khintchine}
and treating the series as a Riemann approximation of an integral we see that
\begin{equation}\label{rhok-limit}
\lim_{t\downarrow0}\log \BE\,\ee^{-ut(\rho_k(t)-\BE\rho_k(t))}=-\int_0^\infty (1-\ee^{-ux}+ux)\frac{\ee^{-kx}}{x}dx
\end{equation}
which is a L\'evy--Khintchine type representation of the centered exponential distribution with mean $1/k$, reflecting the fact that
$t\rho_k(t)$ converges in distribution to $\eps_{1,1}/k$ as $t\downarrow0$.   

Summation of \eqref{rhok-limit} over odd and even $k$ gives
\begin{align}\label{ello-limit}
\lim_{t\downarrow0}\log \BE\,\ee^{-ut(\ell_o(\la(t))-\BE\ell_o(\la(t)))}=-\int_0^\infty (1-\ee^{-ux}+ux)\frac{\ee^{-x}}{x(1-\ee^{-2x})}dx,\\
\label{elle-limit}
\lim_{t\downarrow0}\log \BE\,\ee^{-ut(\ell_e(\la(t))-\BE\ell_e(\la(t)))}=-\int_0^\infty (1-\ee^{-ux}+ux)\frac{\ee^{-2x}}{x(1-\ee^{-2x})}dx.
\end{align}
Interchange of limit, summation and integration is justified by the dominated convergence. It remains to find a sufficiently fine asymptotics 
for $\BE\ell_o(\la(t))$ and $\BE\ell_e(\la(t))$ as $t\downarrow0$ and to identify the limiting distributions described in \eqref{ello-limit},~\eqref{elle-limit}
by their L\'evy measures to prove a limit theorem in the grand canonical ensemble.

The fist task can be accomplished in the simplest way by appealing to results of \cite{FGD}.  We apply Theorem 5~\cite{FGD}
for the base function $g(x)=\frac{\ee^{-x}}{1-\ee^{-x}}$ with the Mellin transform $g^*(v)=\zeta(v)\Gamma(v)$ and the Dirichlet series 
\begin{equation}\label{Dirichlet}
D_o(v)=\sum_{k\text{ odd}}k^{-v}=\left(1-2^{-v}\right) \zeta (v)\qquad\text{and}\qquad D_e(v)=\sum_{k\text{ even}}k^{-v}=2^{-v} \zeta (v)
\end{equation} 
correspondingly, where $\zeta(v)$ is the Riemann zeta function. It gives, with $\gamma$ denoting the Euler--Mascheroni constant, that 
\begin{align}
\BE\,\ell_o(\la(t))=\sum_{k\text{ odd}}\frac{\ee^{-kt}}{1-\ee^{-kt}}=\frac{1}{2t}\log\frac{1}{t}+\frac{\gamma+\log2}{2t}+O(1),\\
\BE\,\ell_e(\la(t))=\sum_{k\text{ even}}\frac{\ee^{-kt}}{1-\ee^{-kt}}=\frac{1}{2t}\log\frac{1}{t}+\frac{\gamma-\log2}{2t}+O(1)
\end{align}
as $t\downarrow0$, because one has
\begin{align*}
D_o(v)g^*(v)&=\frac{1}{2 (v-1)^2}+\frac{\gamma +\log (2)}{2 (v-1)}+O(1),\quad v\to 1,\\
D_e(v)g^*(v)&=\frac{1}{2 (v-1)^2}+\frac{\gamma -\log (2)}{2 (v-1)}+O(1),\quad v\to 1,
\end{align*}
and the only other singular point of $D_o(v)g^*(v)$ and $D_e(v)g^*(v)$ with $\Re v\ge0$
is a simple pole of $D_e(v)g^*(v)$ at $v=0$.

Let $L_o$ and $L_e$ be two random variables with the Laplace transforms as the limits 
\eqref{ello-limit},~\eqref{elle-limit}:
\begin{align}\label{Lo}
\BE\,\ee^{-uL_o}&=\exp\Bigl(-\int_0^\infty (1-\ee^{-ux}+ux)\frac{\ee^{-x}}{x(1-\ee^{-2x})}dx\Bigr)
=\frac{2^u e^{\frac{\gamma  u}{2}} \Gamma \left(\frac{1+u}{2}\right)}{\sqrt{\pi }},\\
\label{Le}
\BE\,\ee^{-uL_e}&=\exp\Bigl(-\int_0^\infty (1-\ee^{-ux}+ux)\frac{\ee^{-2x}}{x(1-\ee^{-2x})}dx\Bigr)=e^{\frac{\gamma  u}{2}} \Gamma \left(1+\tfrac{u}{2}\right).
\end{align}
The closed form expressions can be obtained using Malmst\'en's 
integral formula for $\log\Gamma(z)$, see, e.g., \cite[Ex.~12.3.3]{WW}; we
omit the details. The expression \eqref{Le} is the Laplace transform of the
Gumbel distribution with location parameter $-\gamma/2$ and scale $\frac12$:  
\[
\BE\,\ee^{-uL_e}=e^{\frac{\gamma  u}{2}} \Gamma \left(1+\tfrac{u}{2}\right)=2\int_{-\infty}^\infty 
\ee^{-ux-2x-\gamma -\ee^{-2x-\gamma}}dx,
\]
which can be readily seen with the change of variable $y=\ee^{-2x-\gamma}$. Hence
the cumulative distribution function of $L_e$ is $H_e(x):=\ee^{-\ee^{-2x-\gamma}}$,
$x\in\BR$. Similarly,
\[
\BE\,\ee^{-uL_o}=
\frac{2^u e^{\frac{\gamma  u}{2}} \Gamma \left(\frac{1+u}{2}\right)}{\sqrt{\pi }}=\frac{2}{\sqrt{\pi}}\int_{-\infty}^\infty 
 \ee^{-u x - x-\gamma/2-\log 2 - \ee^{-2 x-\gamma-2\log 2}}
dx, 
\]
that is $L_o$ has the (apparently unnamed) distribution with 
the cumulative distribution function $H_o(x):=\operatorname{Erfc}\bigl(\ee^{-x-\gamma/2}/2)$, $x\in\BR$.
Hence, combining the above results we obtain
\begin{equation}\label{elloe-tlim}
\lim_{t\downarrow0}
\BP\{t\ell_o(\la(t))-\tfrac{1}{2}\log\tfrac{1}{t}-\tfrac{\gamma+\log2}{2}\le x,
t\ell_e(\la(t))-\tfrac{1}{2}\log\tfrac{1}{t}-\tfrac{\gamma-\log2}{2}\le y\}
=H_o(x)H_e(y).
\end{equation}

Passing in \eqref{elloe-tlim} to the time $s=F(t)$ defined in \eqref{F(t)}
as described in Section~\ref{sec:Ppp}, we see that the assumptions
of Theorem~\ref{th:equivalence} hold with $X_1(\la)=\ell_o(\la)$, 
$X_2(\la)=\ell_e(\la)$ and
\[
\ta_1(s)=\ta_2(s)=\frac{3s}{\pi^2}\log\frac{6s}{\pi^2},
\qquad\tb_1(s)=\tb_2(s)=\frac{\pi^2}{6s}.
\]
Condition~\eqref{a-slow} is easy to check for such $\ta_i,\tb_i$ taking
$u(s)=s^{1/8}$, say, and the conclusion of Theorem~\ref{th:equivalence}
is translated into \eqref{elloe-joint}, \eqref{alpha} by simple but tedious 
manipulations with shifts.  
\end{proof}

\begin{remark*}
In terms of the Laplace transforms, the fact that $\ell(\la)=\ell_o(\la)+\ell_e(\la)$
is equivalent to the Legendre duplication formula $\Gamma(2z)=2^{2z-1}\Gamma(z)\Gamma(z+1/2)/\sqrt{\pi}$.  Applied with $z=(1+u)/2$
and multiplied by $\ee^{\gamma u}$,
the right-hand side becomes the limiting Laplace transform of $\ell_o(\la)+\ell_e(\la)$
while the left-hand side is the Laplace transform of the Gumbel distribution \eqref{EL-Gumbel}.
\end{remark*}

\section{Concluding remarks}\label{sec:remarks}
The partition growth process presented in this paper seems to be 
a new approach to study of random partitions. The explicit construction
using independent exponentials described in Section~\ref{sec:gce} 
can be generalized by taking another distribution for the double
sequence, or even by taking differently distributed independent random variables.  
This approach can cause some difficulties, because we 
used explicitly the memoryless property, which characterizes the 
exponential distribution. Perhaps this is not very important and
some other argument can be used. However,
the marked Poisson point process approach of
Section~\ref{sec:Ppp} looks more promising. It seems that in such a way one can
construct a growth process 
leading to any 
multiplicative distribution on integer partitions, introduced
by A.~Vershik in~\cite{V-FAA} as a far-reaching generalization of the uniform measure.

\section*{Acknowledgments} 
This work was partially supported by the Russian Science Foundation, grant 21-11-00141.
The author thanks two anonymous referees for careful reading and useful advises which helped to improve this paper.

\appendix
\section{Calculation of the visit probability}\label{sec:appA}

In this Appendix we aim to prove the last equality in \eqref{g-zero-n}.
We follow a general procedure to decompose a function into infinite series described in \cite[Ch.~7.4]{WW}.
Consider a function
\[
g(z)=\frac{\pi\,\sinh(\frac{\pi z}{6})}{\sqrt{3}(2\cosh(\frac{\pi z}{3})-1)}.
\]
It is a meromorphic function in $\BC$, with an infinite sequence of simple poles 
at purely imaginary points $z=(6k\pm 1)\ii$, $k\in\BZ$, due to zeroes of $2\cosh(\tfrac{\pi z}{3})-1$ at these points.  
It is easy to see that 
\begin{equation}\label{g-res}
\zeta_k\eqdef\operatorname*{Res}\limits_{z=(6k+1)\ii}g(z)=\operatorname*{Res}\limits_{z=(6k-1)\ii}g(z)=\frac{(-1)^{k}}{2},\qquad k\in \BZ.
\end{equation}

Define a sequence of concentric circles $\CC_m=\{z\in\BC:|z|=6m+3\}$ for $m=1,2,\dots$.  
Then $|g(z)|$ attains its maximum on $\CC_m$ for purely imaginary $z$:
\[
\max_{z\in \CC_m}\bigl|g(z)\bigr| = \bigl|g(\pm(6m+3)\ii)\bigr|=\frac{\pi}{3\sqrt{3}}.
\]
Consider a point $x$ which is not a pole of $g(z)$.  Then $g(z)/(z-x)$ has a simple pole at $z=x$ 
with residue $g(x)$ and a sequence of poles $(6k\pm1)\ii$ of $g(z)$ with residues $\zeta_k/((6k\pm 1)\ii-x)$,
$k\in\BZ$.
Hence, if $6m+3>|x|$, by the Cauchy residue theorem we obtain
\begin{align*}
\frac{1}{2\pi\ii}\oint_{\CC_m}\frac {g(z)}{z-x}dz 
&=g(x)+\sum_{|k|\le m}\frac{(-1)^{k}}{2((6k+ 1)\ii-x)}+\sum_{|k|\le m}\frac{(-1)^{k}}{2((6k- 1)\ii-x)}\\
&=g(x)+\sum_{|k|\le m}\frac{(-1)^{k+1}x}{(6k+ 1)^2+x^2}\,.
\end{align*}
To obtain the last line we have changed $k$ to $-k$ in the second sum and joined two sums together.
On the other hand,
\begin{equation}\label{two-integrals}
\frac{1}{2\pi\ii}\oint_{\CC_m}\frac {g(z)}{z-x}dz =
\frac{1}{2\pi\ii}\oint_{\CC_m}\frac {g(z)}{z}dz +\frac{1}{2\pi\ii}\oint_{\CC_m}\frac {xg(z)}{z(z-x)}dz .
\end{equation}
Since $g(z)$ is bounded on $\CC_m$, the last integral vanishes as $m\to\infty$ along natural numbers.  
The first integral on the right-hand side of \eqref{two-integrals} can be evaluated again by the Cauchy 
residue theorem: since $g(0)=0$, 
\[
\frac{1}{2\pi\ii}\oint_{\CC_m}\frac {g(z)}{z}dz=\sum_{|k|\le m}\frac{\zeta_k}{(6k+1)\ii} +\sum_{|k|\le m}\frac{\zeta_k}{(6k-1)\ii}=0
\]
by using $\zeta_k=\zeta_{-k}$ and changing $k$ to $-k$ in the second sum.  Therefore we obtain 
\begin{equation}\label{g(x)}
g(x)+\sum_{k\in \BZ}\frac{(-1)^{k+1}x}{(6k+ 1)^2+x^2}=0.
\end{equation}
It remains to multiply \eqref{g(x)} by $24n/x$,  take $x=\sqrt{24n-1}$ and notice that
\[
\frac{1}{n+\frac{k(3k+1)}{2}}
=\frac{24}{(6k+1)^2+24n-1}
\]
to obtain \eqref{g-zero-n}.

\section{List of notation used}\label{sec:appB}

\begin{center}
\begin{tabular}{p{0.15\textwidth} p{0.65\textwidth} c}
Notation & Description & Reference \\
\hline 
$C_k(\la)$\vphantom{$\hLa$} & Number of parts $k$ in a partition $\la$ & \eqref{N-via-Rk} \\
$\eps_{k,j}$ & Independent standard exponential random variables & p.~\pageref{eps} \\
$F(t)$ & Tail of intensity measure of the Poisson point process $\Pi$ & \eqref{F(t)} \\
$F_k(t)$ & Tail of intensity measure of the Poisson point process $\Pi_k$ & \eqref{Fk(t)} \\
$f(t)$ & Density of the Poisson point process $\Pi$; $f(t)=-F'(t)$ & \eqref{f-def} \\
$f_k(t)$ & Density of the Poisson point process $\Pi_k$; $f_k(t)=-F_k'(t)$ & \eqref{fk-def} \\
$g^\tau(\la)$, $g^\tau(n)$ & The probability that the chain $(\la(t))_{t\ge\tau}$ visits
$\la\in\CP_n$& \eqref{g-la}, \eqref{gn}\\
$g(\la)$, $g(n)$ & The probability that the chain $(\la(t))_{t>0}$ visits
$\la\in\CP_n$& \eqref{g-zero-n}, \eqref{gn}\\
$\gamma(n)$ & The probability that $(\la(t))_{t>0}$ visits
some partition in $\CP_n$ & \eqref{visit-prob}\\
$J(\tau)$ & The number of jumps of $\bigl(\la(t)\bigr)_{t\ge\tau}$ & p.~\pageref{sec:jumps} \\
$\ell(\la)$ & Length of a partition $\la$, that is number of its parts & \eqref{N-via-Rk} \\
$\bigl(\la(t)\bigr)_{t\ge\tau}$ & Decreasing (forward)  partition process & p.~\pageref{la(t)}\\
$\bigl(\tla(s)\bigr)_{s\ge 0}$ & Increasing partition process in time $s=F(t)$ with $F$ defined in \eqref{F(t)}: $\tla(s)=\la(F^{-1}(s))$ & p.~\pageref{tla(s)}\\
$(\Lambda^\tau_j)_{j=0,\dots,J(\tau)}$ & Decreasing (forward) jump chain of $\bigl(\la(t)\bigr)_{t\ge\tau}$ & p.~\pageref{sec:jumps} \\
$(\Lambda_j)_{j\le 0}$ & Decreasing (forward) jump chain indexed by non-positive integers with terminal state $\Lambda_0=\varnothing$
and no entry state & p.~\pageref{Lambda-j} \\
$(\hLa^\tau_j)_{j=0,\dots,J(\tau)}$ & Increasing (backward) jump chain of $\bigl(\la(t)\bigr)_{t\ge\tau}$ & \eqref{hat-Lambda} \\
$(\hLa_j)_{j\ge0}$ & Increasing (backward) jump chain of $\bigl(\la(t)\bigr)_{t>0}$ & Th.~\ref{th:hat-Lambda} \\
$M_n^{\text{uni}}$ & The uniform measure on the set $\CP_n$ & p.~\pageref{Mnuni} \\
$M_q^{\text{gc}}$ & Fristedt's measure on $\CP$ with parameter $q\in(0,1)$ & \eqref{Mgc}\\
$N(\la)$ & The weight of a partition $\la$, that is $\sum_j\la_j$ & \eqref{N-via-Rk}\\
$\CP$ & The set of all integer partitions, $\CP=\cup_{n\ge 0}\CP_n$ & p.~\pageref{CP} \\
$\CP_n$ & The set of integer partitions of $n$ & p.~\pageref{CP} \\ 
$p_{\la,\mu}$ & Transition probabilities of the decreasing (forward) jump chains $(\Lambda^\tau_j)$ and $(\Lambda_j)$ & Prop.~\ref{prop:jump-forward}\\ 
$\hp_{\mu,\la}$ & Transition probabilities of the increasing (backward) jump chain $(\hLa^\tau_j)$ and $(\hLa_j)$ & \eqref{edge-probability}\\ 
$\Pi$ & The Poisson point process of $\la(t)$ jump times &  p.~\pageref{Pi-def}\\
$\Pi_k$ & The Poisson point process of $\rho_k(t)$ jump times &  p.~\pageref{Pik-def}\\
$\tPi$ & The standard Poisson point process of $\tla(s)$ jump times &  p.~\pageref{tPi-def}\\
$q_{n,m}$ & Transition probabilities of decreasing Markov chain $N(\Lambda_j)$ & \eqref{qnm}\\
$\hq_{n,m}$ & Transition probabilities of increasing Markov chain $N(\hLa_j)$ & \eqref{hq}\\
$\rho_k(t)$ & Counting process for parts $k$ in $\la(t)$ & \eqref{Rk-def}\\ 
$\trho_k(s)$ & Counting process for parts $k$ in $\tla(s)$; $\trho_k(s)=\rho_k(F^{-1}(s))$ & \eqref{trhok-def}\\ 
\end{tabular}
\begin{tabular}{p{0.15\textwidth} p{0.65\textwidth} c}
$\sigma_1(n)$ & Sum of divisors function & \eqref{hq} \\
$Y_\la$ & The Young diagram of a partition $\la$ & \eqref{def-Y} \\
$\varnothing$ & The empty partition of $0$ or the empty set, depending on the context 
& p.~\pageref{empty-partition} \\
$\nearrow$, $\nwarrow$ & $\mu\nearrow\la$ or $\la\nwarrow\mu$ means that $Y_\la$ is obtained from $Y_\mu$ by adding one box & p.~\pageref{nearrow}\\
$\vdash$ & $\la\vdash n$ means that $\la$ is a partition of $n$ & p.~\pageref{vdash}\\
$\prec$, $\preccurlyeq$ & $\mu\prec\la$ means that one can obtain $\la $ from $\mu$ adding further parts to~$\mu$ and not changing or removing existing parts & p.~\pageref{preccurlyeq}
\end{tabular}
\end{center}


\begin{thebibliography}{99}

\bibitem{Aigner} Aigner, M. \textit{A Course in
Enumeration}. Springer-Verlag, Berlin, Heidelberg, 2007

\bibitem{Andrews} Andrews, G.\;E. \textit{The Theory of Partitions}. Addison--Wesley Publishing Company, 1976.

\bibitem{Baez-Duarte}
B\'aez--Duarte, L. Hardy--Ramanujan's asymptotic formula for partitions and the central limit theorem. \textit{Adv. Math.} \textbf{125}, no.~1, 
114--120 (1997).


\bibitem{Bogachev}
Bogachev, L.\;V.
Unified derivation of the limit shape for multiplicative ensembles of random integer partitions with equiweighted parts.
\textit{Random Structures Algorithms} \textbf{47}, no.~2, 227--266 (2015).


\bibitem{BogYak}
Bogachev, L.\;V., Yakubovich, Y.\;V. Limit shape of minimal difference partitions and fractional statistics. 
\textit{Commun. Math. Phys.} \textbf{373}, 1085--1131 (2020). 


\bibitem{Corteel}
Corteel, S., Pittel, B., Savage, C.\;D.,  Wilf, H.\;S. On the multiplicity of parts in a random partition, 
\textit{Random Structures Algorithms} \textbf{14}, no.~2, 185--197 (1999).

\bibitem{DeSalvo}
DeSalvo, S. Will the real Hardy--Ramanujan formula please stand up? \textit{Integers} \textbf{21}, \#A88 (2021).

\bibitem{EL}
Erd\H{o}s, P., Lehner,  J. The distribution of the number of summands in the partitions of a positive integer,
\textit{Duke Math. J.}, \textbf{8}, 335--345 (1941).
%
%



\bibitem{Eriksson}
Eriksson, K., Sj\"ostrand, J.
Limiting shapes of birth-and-death processes on Young diagrams,
\textit{Adv. Appl. Math.} 
\textbf{48}, no.~4,
575--602 (2012).

\bibitem{Fatkullin}
 Fatkullin, I., Sethuraman, S.,  Xue J. On hydrodynamic limits of Young diagrams. \textit{Electron. J. Probab.} \textbf{25} 1--44 (2020). 

\bibitem{FGD}
Flajolet, P., Gourdon, X., Dumas, P. 
Mellin transforms and asymptotics: harmonic sums. 
\textit{Theoret. Comput. Sci.} \textbf{144}, no. 1--2, 3--58 (1995).


\bibitem{Fristedt} Fristedt, B. The structure of random partitions of large integers. \textit{Trans. AMS} \textbf{337}, no.~2, 703-–735 (1993).


\bibitem{Gnedin-Yak}
Gnedin, A., Yakubovich, Y. On the number of collisions in  $\Lambda$-coalescents, \textit{Electron. J. Probab.} \textbf{12}), 1547--1567, (2007)


\bibitem{Goh-Hitczenko} Goh, W.\;M.\;Y., Hitczenko, P. Random partitions with restricted
part sizes. \textit{Random Structures Algorithms} \textbf{32}, 440--462 (2008).


\bibitem{Granovsky}
Granovsky, B., Stark, D., Erlihson, M. Meinardus’ theorem on weighted partitions: Extensions and a probabilistic proof, \textit{Adv.
Appl. Math.} \textbf{41}, 307--328 (2008).

\bibitem{Hardy-Ramanujan}
Hardy, G.\;H., Ramanujan, S. Asymptotic formul\ae in combinatory analysis. \textit{Proc. Lond.
Math. Soc. (3)}, \textbf{2}, no. 1, 75--115  (1918).

\bibitem{Hardy-Wright}
Hardy, G.\;H., Wright, E.\;M. \textit{An Introduction to the Theory of Numbers}, Fourth ed., Clarendon Press, Oxford, 1975.

\bibitem{Hunt}  Hunt, G.\;A.
Markov chains and Martin boundaries,
\textit{Illinois J. Math.} \textbf{4},
313--340  (1960).

\bibitem{James}  James, G.\;D. \textit{The Representation Theory of Symmetric Groups},
Lecture Notes in Math., Vol. 682, Springer-Verlag, New York, NY,
1978.

\bibitem{Kerov93} Kerov, S.\;V. Transition Probabilities of Continual Young Diagrams and Markov Moment
Problem, \textit{Funktsion. Anal. i Prilozhen.} \textbf{27}, no. 2, 32--49 (1993); English transl.:
\textit{Funct. Anal. Appl.} \textbf{27}, 104--117 (1993). 

\bibitem{Kerov99} Kerov, S.\;V. 
A differential model for the growth  of Young diagrams. Proceedings of the St. Petersburg Mathematical Society, Vol. IV, 111--130, Amer. Math. Soc. Transl. Ser. 2, 188, Amer. Math. Soc., Providence, RI, 1999. 

\bibitem{Kingman} Kingman, J.\;F.\;C. \textit{Poisson Processes}. Clarendon Press, 1993.

\bibitem{Krapivsky}
 Krapivsky, P.\;L. Stochastic dynamics of growing Young diagrams and their limit shapes. \textit{J. Stat. Mech.} 013206 (2021).

\bibitem{Lovasz-Plummer}
Lov\'{a}sz, L.,  Plummer, M.\;D. \textit{Matching Theory.} North-Holland Mathematics Studies. Elsevier Science, 1986.

\bibitem{Mutafchiev}
Mutafchiev, L.\;R.
On the size of the Durfee square of a random integer partition,
\textit{J. Comput. Appl. Math.} \textbf{142}, no.~1,
173--184 (2002).

\bibitem{Norris} Norris, J.\;R. \textit{Markov Chains}. Cambridge Univ. Press, 1997.

\bibitem{Olsh-Ivanov} Ivanov, V., Olshanski, G. Kerov’s central limit theorem for the
Plancherel measure on Young diagrams. In S. Fomin (ed.) \textit{Symmetric Functions 2001: Surveys of Developments and Perspectives}
(NATO Science Series II. Mathematics, Physics and Chemistry. Vol. 74), Kluwer, 2002, pp. 93--151.


\bibitem{Petrov} Petrov, F. Limits shapes of Young diagrams. Two elementary approaches.
\textit{Zapiski Nauchn. Sem. POMI} \textbf{370}, 111--131 (2009) (in Russian); English transl. in \textit{J. Math. Sci.} \textbf{166}, no.~1, 63--74 (2010).

\bibitem{CSP} Pitman, J. \textit{Combinatorial Stochastic Processes}. \'Ecole d’\'et\'e de Probabilit\'es de Saint-
Flour XXXII – 2002, Lecture Notes in Math., vol. 1875, Springer, 2006.

\bibitem{Pittel}
Pittel, B. On a likely shape of the random Ferrers diagram, \textit{Adv. Appl. Math.} 18, no.~4, 432--488 (1997).

%

\bibitem{Stam} Stam, A.\;I. Independent Poisson processes generated by record values and inter-record
times. \textit{Stochastic Process. Appl.}, \textbf{19}, no. 2, 315--325  (1985).

\bibitem{Strassen}
Strassen, V. The existence of probability measures with given
marginals. \textit{Ann. Math. Stat.} \textbf{36}, 423--439  (1965).

\bibitem{Su15}
Su, Z. Normal convergence for random partitions with multiplicative measures.
\textit{Th. Probab.Appl.} \textbf{59}, no.~1, 40--69  (2015).

\bibitem{Temperley}  Temperley, H.\;N.\;V. Statistical mechanics and the partition of numbers. II. The form of crystal surfaces. \textit{Proc. Cambridge
Philos. Soc.} \textbf{48}, 683--697 (1952). 
 

\bibitem{Vershik-Kerov77}
Vershik, A.\;M., Kerov, S.\;V. Asymptotics of the Plancherel measure of the symmetric group and the limiting form of Young tableaux.
\textit{Dokl. AN SSSR}, \textbf{233}, no.~6, 1024--1027 (1977); English transl. in \textit{Soviet Math. Dokl.} \textbf{18}, 527--531 (1977).

\bibitem{V-FAA}
 Vershik, A.\;M. Statistical mechanics of combinatorial partitions, and their limit configurations (Russian), \textit{Funktsional. Anal. i
Prilozhen.} \textbf{30}, no. 2, 19--39 (1996);  English transl. in \textit{Funct. Anal. Appl.} \textbf{30}, no. 2, 90--105 (1996). 

\bibitem{V-UMN}
 Vershik, A.\;M. Limit distribution of the energy of a quantum ideal gas from the point of view of the theory of partitions of natural numbers. 
 \textit{Russian Math. Surveys} \textbf{52}, no. 2, 379--386 (1997).

\bibitem{WW}
Whittaker, E.\;T., Watson, G.\;N. \textit{A Course of Modern Analysis}, 4th ed., Cambridge Univ. Press, 1958.
 
\bibitem{Yakubovich12} Yakubovich, Y.  
  Ergodicity of multiplicative statistics. \textit{J. Comb. Theory Ser. A} \textbf{119}, no. 6, 1250--1279 (2012). 
 
\bibitem{Yakubovich23} 
Yakubovich, Yu. V. Moments of random integer partitions.
\textit{Zapiski Nauchn. Sem. POMI} \textbf{525}, 161--183 (2023) (in Russian).
English transl. to appear in \textit{J. Math. Sci.}
 
\end{thebibliography}
\end{document}